\documentclass[11pt,paper=letter]{scrartcl}

\usepackage{graphicx}
\usepackage{amsmath}
\usepackage{amsfonts}
\usepackage{amssymb}
\usepackage{bm}

\usepackage{times}
\usepackage{xcolor}

\usepackage[numbers,sort]{natbib}
\usepackage{xspace}
\usepackage{xcolor}
\usepackage{amssymb}
\usepackage{enumerate}
\usepackage{verbatim}
\usepackage{subfigure}
\usepackage{bbm}
\usepackage{xr}
\usepackage{threeparttable}

\usepackage[colorlinks,linkcolor=blue,citecolor=red]{hyperref}
\usepackage{hypernat}
\usepackage[bf,format=plain,singlelinecheck=false]{caption}

\usepackage{algorithm}
\usepackage{algorithmic}
\usepackage{setspace,url,bm,amsmath} 

\def\T{{ \mathrm{\scriptscriptstyle T} }}
\def\v{{\varepsilon}}
\def\e{{\eta}}

\def\var{\textnormal{var}}
\def\ni{n_i}
\def\sumij{\sum_{i=1}^{B}\sum_{j=1}^{\ni}}
\def\sumi{\sum_{i=1}^{B}}
\def\sumj{\sum_{j=1}^{\ni}}

\def\maxi{\max\limits_{i=1,\ldots,B}}

\def\maxj{\max\limits_{j=1,\ldots,\ni}}
\def\hattaui{\hat \tau_{i \cdot } }
\def\taui{ \tau_{i\cdot} }
\def\nti{n_{1i}}
\def\nci{n_{0i}}
\def\unadj{\textnormal{unadj}}
\def\ols{\textnormal{ols}}
\def\olsint{\textnormal{ols\_int}}
\def\tauunadj{\hat \tau_{\unadj}}
\def\tauols{\hat \tau_{\ols}}
\def\tauiols{\hat \tau_{i,\ols}}

\def\tauiolsint{\hat \tau_{i,\olsint}}
\def\tauolsint{\hat \tau_{\olsint}}

\def\S{S}
\def\s{s}
\def\X{X}
\def\x{X}

\newcommand\revise[1]{\textcolor{black}{#1}}

\def\Sunadj{\sigma_{\unadj}}
\def\sunadj{\hat \sigma_{\unadj}}

\def\argmin{\mathop{\arg\min}}

\newtheorem{theorem}{Theorem}
\newtheorem{lemma}{Lemma}
\newtheorem{corollary}{Corollary}
\newtheorem{condition}{Condition}
\newtheorem{remark}{Remark}
\newtheorem{proposition}{Proposition}

\newtheorem{proof}{Proof}

\title{Regression-adjusted average treatment effect estimates in stratified randomized experiments }

\author{
	Hanzhong Liu\thanks{Center for Statistical Science and Department of Industrial Engineering, Tsinghua University, Beijing, 100084, China. Email:lhz2016@tsinghua.edu.cn}
	\and
	Yuehan Yang\thanks{School of Statistics and Mathematics, Central University of Finance and Economics, Beijing, 102206, China. Email:yyh@cufe.edu.cn}
}

\begin{document}

\maketitle


\begin{abstract}
  \noindent	\footnotesize\textbf{Abstract}

  \noindent	 Researchers often use linear regression to analyse randomized experiments to improve treatment effect estimation by adjusting for imbalances of covariates in the treatment and control groups. Our work offers a randomization-based inference framework for regression adjustment in stratified randomized experiments. Under mild conditions, we re-establish the finite population central limit theorem for a stratified experiment. We prove that both the stratified difference-in-means and the regression-adjusted average treatment effect estimators are consistent and asymptotically normal. The asymptotic variance of the latter is no greater and is typically lesser than that of the former. We also provide conservative variance estimators to construct large-sample confidence intervals for the average treatment effect.

  \ \\

  \noindent \textbf{Keywords}

  \noindent Blocking; Randomized-block design; Randomized experiments; Randomization-based inference; Stratified sampling.
\end{abstract}


\section{Introduction}

Randomization experiments are the gold standard for drawing causal inference. The Neyman--Rubin potential outcomes model \citep{Neyman:1923,Rubin:1974} is widely used to analyse the experimental results. Under the stable unit treatment value assumption \revise{\citep{Rubin:1980}}, the average treatment effect can be estimated without bias by the simple difference-in-means estimator. In practice, randomization is often stratified according to some categorical variables. Appropriate stratification can better balance baseline covariates and improve the estimation efficiency \citep{fisher1926,imai2008variance,miratrix2013,imbens2015causal}.

Apart from stratification variables, researchers often observe many other baseline covariates. Even when carefully designed, the probability of some baseline covariates exhibiting imbalance between treatment and control groups is high \citep{fisher1926, morgan2012rerandomization}. Regression adjustment is a common strategy to find a more efficient estimator. In completely randomized experiments, \citet{lin2013} showed that regression adjustment can yield a consistent and asymptotically normal estimator whose asymptotic variance is no greater than that of the difference-in-means estimator. \revise{\citet{bloniarz2015lasso} used the lasso to perform regression adjustment for a high-dimensional setting.} \citet{fogarty2018finely,fogarty2018}  discussed the advantages of regression adjustment in two special stratified randomized experiments, finely stratified and paired experiments. The above results \revise{were} obtained under a randomization-based inference framework \citep{kempthorne1955randomization, fcltxlpd2016, zhao2016randomization}. This framework assumes the correctness of neither a fitted regression model nor a super-population model. Instead, the potential outcomes and covariates are assumed to be fixed quantities and the randomness comes only from the treatment assignment. The experimental units consist of a finite population, which is embedded into a hypothetical infinite sequence of finite populations with increasing sizes.

\revise{Recently, stratified randomized experiments were discussed under the Neyman--Rubin model and randomization-based inference framework  \citep{Imai2008, imbens2015causal, Higgins2015, Schochet2016, Pashley2017}, however no formal asymptotic results have been established. When the number of strata  is fixed with their sizes going to infinity, the asymptotic results of \citet{fcltxlpd2016}  can be directly extended to stratified randomized experiments. However, it is not clear whether or not their results hold when the number of strata tends to infinity. We first fill in these gaps by showing that the simple stratified difference-in-means estimator is consistent and asymptotically normal in the asymptotic regime where the total number of experimental units $N$ tends to infinity, allowing the number of strata or their sizes to increase with $N$. Second, we discuss two regression-adjusted average treatment effect estimators to improve the estimation efficiency. Under mild conditions, we show that they are also consistent and asymptotically normal. For the first estimator, we consider a general asymptotic regime, allowing the number of strata to tend to infinity. It extends the asymptotic results of regression adjustment from finely stratified and paired experiments \citep{fogarty2018finely,fogarty2018} to experiments with broader combinations of the number of strata and their sizes. In addition, we show that the first estimator improves, or at least does not hurt, precision when compared with the stratified difference-in-means estimator, provided the proportion of treated units is asymptotically the same across strata. The second estimator, whose asymptotic normality is a straightforward extension of \citet{fcltxlpd2016}, is more efficient than the first one when  the number of strata is fixed with sizes going to infinity. Finally, we provide conservative variance estimators to construct large-sample conservative confidence intervals for the average treatment effect. }

\section{Stratified randomized experiments}

Consider a stratified randomized experiment with $N$ units consisting of $B$ strata. Within the strata, independent completely randomized experiments are conducted, allowing different relative sizes of treatment and control groups. The stratum $i$ contains $\ni$ units, $n_i \geq 2$, such that $\sum_{i=1}^{B} \ni = N$, of whom $\nti$ are sampled without replacement and receive the active treatment, while the remaining $\nci = \ni - \nti$ receive the control. \revise{Here and in what follows, we use `$1$' to denote treatment and `$0$' to denote control}. Let $Z_{ij}$ be an indicator of whether or not the $j$th unit in stratum $i$ is sampled and receives the treatment, and $\sum_{j=1}^{\ni} Z_{ij} = \nti$. In stratified random sampling, the probability that the indication vector $Z = (Z_{11}, \ldots, Z_{Bn_B} )$ takes a particular value $(z_{11},\ldots,z_{Bn_B})$ is $\prod_{i=1}^{B} \nti! \nci! / \ni !$, where $\sum_{j =1}^{\ni} z_{ij} = \nti$, for $i = 1,\ldots, B$. The total treated sample size is $ \sumi \nti$. For each unit $j$ in stratum $i$, there are two fixed potential outcomes, \revise{$y_{ij}(1)$ and $y_{ij}(0)$, representing the outcomes of this unit receiving and not receiving the treatment, respectively}. Unit or individual level treatment effect, $\tau_{ij}$, is defined as the comparison of these two potential outcomes, \revise{ $\tau_{ij} = y_{ij}(1) - y_{ij} (0)$}.  Since each unit is either in the treatment or in the control group, but not both, potential outcomes \revise{$y_{ij}(1)$ and $y_{ij}(0)$} cannot be observed simultaneously, and hence, $\tau_{ij}$ is not identifiable without strong model assumptions. However, the average across all experimental units is estimable. It is often called the average treatment effect and denoted as 
$$ \tau  = \sumij \tau_{ij} /N = \sum_{i=1}^{N} c_i \tau_{i\cdot}, $$
where $c_i = \ni/N $ is the proportion of stratum size, and  $\tau_{i\cdot} = \sum_{j=1}^{\ni} \tau_{ij}/\ni$ is the  average treatment effect in stratum $i$. The observed outcome is 
$$y_{ij} = Z_{ij} y_{ij}(1) + ( 1 - Z_{ij} ) y_{ij}(0) . $$

For stratum $i$, let \revise{$\hat y_{i\cdot}(1)$ and $\hat y_{i\cdot}(0)$} be the sample mean of  the observed outcomes in the treatment and control groups, respectively. A natural estimator for  $\tau_{i\cdot}$ is the difference of these two sample means,  \revise{$ \hat \tau_{i \cdot} = \hat y_{i\cdot}(1)  - \hat y_{i\cdot}(0)$}. The stratified difference-in-means estimator for the total average treatment effect $\tau$ is 
$$ \tauunadj = \sumi c_i \hat \tau_{i \cdot} ,$$
where the subscript `unadj' indicates that the estimator does not adjust imbalances of covariates.

Apart from stratification variables, the unit $j$ in stratum $i$ has a $K$-dimensional  baseline covariate vector, $\X_{ij} =  (\x_{ij1},\ldots, \x_{ijK})$, which is measured in principle before randomization and thus \revise{is} not affected by the treatment assignment. If we treat the covariates as outcomes, then \revise{$\X_{ij}(1) = \X_{ij}(0) = \X_{ij}$}. Both the potential outcomes and  covariates  are fixed quantities. The randomness comes only from the treatment assignment vector $Z$.

For any fixed quantities $a_{ij}(1)$ and $a_{ij}(0)$, where $i = 1,\ldots,B, \ j = 1,\ldots, \ni $,  we use a dot `$\cdot$' to denote the average along that coordinate. In particular, for stratum $i$, we define the stratum-specific population means as 
$$ a_{i\cdot}(1) = \ni^{-1} \sumj a_{ij}(1), \quad  a_{i\cdot}(0) = \ni^{-1} \sumj a_{ij}(0).$$
We use a hat on the population mean to denote the corresponding stratum-specific sample mean, that is, 
$$ \hat a_{i\cdot}(1) = \nti^{-1} \sumj Z_{ij} a_{ij}(1), \quad \hat a_{i\cdot}(0) = \nci^{-1} \sumj (1 - Z_{ij} ) a_{ij}(0) . $$
The stratum-specific population variances are denoted as $\S^2_{ia}(1)$ and $\S^2_{ia}(0) $, and the stratum-specific sample variances are denoted as $\s^2_{ia}(1)$ and $\s^2_{ia}(0)$. Furthermore, let $S_{i \X \X}$  be the population covariance matrix of covariates in  stratum $i$, and denote $S_{i \X y}(1)$ and $S_{i \X y}(0)$ as the population covariance vectors between the covariates $\X_{ij}$ and the potential outcomes $y_{ij}(1)$ and $y_{ij}(0)$ in stratum $i$, respectively. Let $p_i = \nti / n_i$ be the stratum-specific proportion of treated units. All the above defined quantities depend on $N$, but we omit this dependence for notation simplicity. Recalling that $c_i = \ni / N$, let
\revise{
\begin{equation}
\Sunadj^2 = \sumi c_i^2 \Big\{  \frac{ \S^2_{iy}(1) } { \nti } + \frac{ \S^2_{iy}(0) }{ \nci } -  \frac{ \S^2_{i \tau} }{ \ni }  \Big\}. \nonumber
\end{equation}}
We first introduce the following proposition.

\begin{proposition}\citep{miratrix2013,imbens2015causal}
The mean and variance of $\tauunadj$ are $ E(\tauunadj) = \tau$ and $  \var (\tauunadj) = \Sunadj^2. $
\end{proposition}


The stratum-specific variance of  $\tau_{ij}$, $ \S_{i \tau}^2 $, is generally not estimable because we cannot observe $\tau_{ij}$, which equals zero if and only if the stratum-specific treatment effect is constant. The variance of the stratified difference-in-means estimator, $\Sunadj^2$, can be consistently estimated if and only if $\tau_{ij}$ is constant within each stratum. It does not require $\tau_{ij}$ to be constant between strata, which is different from that in completely randomized experiments. In general, the variance $\Sunadj^2$ can be estimated by a Neyman-type conservative estimator \revise{when $\nti, \nci \geq 2$, that is, }
$$
\sunadj^2 = \sumi c_i^2  \Big\{  \frac{ \s^2_{iy}(1) }{ \nti }  +  \frac{ \s^2_{iy}(0) }{ \nci }   \Big\}.
$$

To conduct statistical inference of $\tau$ based on $\tauunadj$, we need to study the sampling distribution, or asymptotic distribution, of $\tauunadj$ induced by the stratified experiment. Finite population central limit \revise{theorem} plays an essential role for randomization-based asymptotic inference \citep{sen1995hajek, robinson1978asymptotic,fcltxlpd2016}.
In what follows, we re-establish \revise{the asymptotic normality result in \citet{bickel1984asymptotic}}, which is similar to the classical Lindeberg--Feller central limit \revise{theorem} for independent random variables. To obtain this result, we propose a more interpretable condition to replace the Lindeberg--Feller condition, although it is slightly stronger than the latter.

\section{Asymptotic normality of stratified difference-in-means estimator}

\subsection{Finite population central limit \revise{theorem}}


\revise{The treatment group can be seen as a stratified random sample from the population $\Pi_{N}(1) = \{ y_{ij}(1): \ i = 1,\ldots, B; \ j = 1,\ldots, \ni\}$. The finite population central limit theorem for a stratified random sample depends crucially on the maximum weighted squared distance of $y_{ij}(1)$ from its stratum-specific average $y_{i \cdot}(1)$:
$$
m_{1N} = \max\limits_{i=1,\ldots,B} \max\limits_{j=1,\ldots,\ni}   \frac{1}{p_i^2}   \frac{ \big\{ y_{ij}(1) - y_{i\cdot}(1) \big\}^2 }{ \sumi c_i \S^2_{iy}(1)  \nci  / \nti } .
$$
Let $ y_{\cdot \cdot }(1) = \sumi \sumj y_{ij}(1) / N = \sumi c_i y_{i \cdot }(1)$ be the population mean and let $ \hat y_{\cdot \cdot} (1) = \sumi c_i \hat y_{i \cdot}(1) $} be the weighted sample mean. Define the weighted population and sample \revise{variances as
\begin{equation}
\sigma^2_{y}(1) = \sumi c_i^2 \S^2_{iy}(1)  \Big( \frac{1}{ \nti } - \frac{1}{\ni } \Big), \quad  \hat \sigma^2_{y}(1) = \sumi c_i^2 \s^2_{iy}(1)   \Big( \frac{1}{ \nti } - \frac{1}{\ni } \Big) , \nonumber
\end{equation}
respectively. Let $\mathcal{N}(0,1)$ be the standard normal distribution.} We have the following results.
\revise{
\begin{proposition}\citep{cochran1977}
\label{prop:mean_sd}
For stratum $i$, the mean and variance of $\hat y_{i\cdot}(1) $ are $E\{ \hat y_{i\cdot}(1) \} = y_{i \cdot}(1)$ and $\var \{ \hat y_{i\cdot}(1) \} = \S^2_{iy}(1)  ( 1/\nti - 1/\ni )$. Furthermore, the mean and variance of $\hat y_{\cdot \cdot} (1) $ are  $E\{ \hat y_{\cdot \cdot} (1)  \} = y_{\cdot \cdot}(1)$ and $\var \{ \hat y_{\cdot \cdot} (1)  \} = \sigma^2_{y}(1)$.
\end{proposition}
}
\begin{theorem}
\label{thm:fclt}
If $m_{1N} /N \rightarrow 0$, then \revise{$ \big\{ \hat y_{\cdot \cdot} (1) - y_{\cdot \cdot }(1) \big\} / \sigma_{y}(1)$  converges in distribution to $\mathcal{N}(0,1)$}, as $N \rightarrow \infty$. Furthermore, if for each stratum $i$,  $2 \leq \nti \leq \ni - 2$, then, \revise{$\hat \sigma^2_{y}(1) / \sigma^2_{y}(1)$} converges to one in probability.
\end{theorem}

\revise{
Theorem~\ref{thm:fclt} suggests a strategy to construct a large-sample confidence interval for the population mean, $y_{\cdot \cdot}(1)$, based on the Normal approximation. When there is only one stratum, $B=1$, the condition $m_{1N} /N \rightarrow 0$ is reduced to that in \citet{fcltxlpd2016} proposed to obtain the finite population central limit \revise{theorem} for a simple random sample. This theorem substantially generalizes the existing result because it allows $B \rightarrow \infty$ as $N \rightarrow \infty$ with fixed $\ni$, although it is a direct extension of \citet{fcltxlpd2016} when $B$ is bounded and $n_i \rightarrow \infty$.
}

\subsection{Asymptotic normality of stratified difference-in-means estimator}

We apply Theorem~\ref{thm:fclt} to derive the asymptotic normality of  stratified difference-in-means estimator $\tauunadj$ for inferring the average treatment effect $\tau$. To present the theoretical \revise{result}, we propose the following conditions.


{\color{black}
\begin{condition}
\label{cond:prop}
There exists constants $p_{i,\infty}$ and $C \in (0, 0.5)$, such that $C < \min_{i= 1,\ldots, B} p_{i, \infty} \leq  \max_{i= 1,\ldots, B} p_{i, \infty} < 1 - C$ and  $\max_{i=1,\ldots,B} | p_i  - p_{i,\infty} |  \rightarrow 0$.
\end{condition}
\begin{condition}
\label{cond:max_dist_y}
The maximum stratum-specific squared distances divided by $N$ tend to zero,
\begin{equation}
\frac{1}{N} \maxi \maxj \big\{ y_{ij}(1) - y_{i\cdot}(1) \big\}^2  \rightarrow 0, \quad  \frac{1}{N} \maxi \maxj \big\{ y_{ij}(0) - y_{i\cdot}(0) \big\}^2  \rightarrow 0. \nonumber
\end{equation}
\end{condition}
\begin{condition}
\label{cond:moments_outcomes}
The weighted variances $\sumi c_i  \S^2_{iy}(1)  / p_i $, $ \sumi c_i  \S^2_{iy}(0)  / (1 - p_i ) $ and $\sumi c_i   \S_{i \tau}^2 $ tend to finite limits, positive for the first two, and the limit of $\sumi c_i \big\{ \S^2_{iy}(1)  / p_i  + \S^2_{iy}(0)  / (1 - p_i ) - \S^2_{ i \tau } \big\}$ is positive.
\end{condition}}
\begin{remark}
\revise{The proportion of treated units $p_i$ is sometimes the same for all strata. In this case, the stratification is beneficial if based on the covariates that are strongly correlated with the outcomes. Otherwise, the stratification may hurt precision. We refer the reader to  \citet{neyman1935statistical}, \citet{miratrix2013}, and \citet{Pashley2017} for more detailed discussion.} Condition~\ref{cond:max_dist_y} is slightly stronger, but simpler and more interpretable, than the Lindeberg--Feller condition proposed by \citet{bickel1984asymptotic}. It is weaker than the fourth moment condition proposed by \citet{fogarty2018finely, fogarty2018}. In Condition~\ref{cond:moments_outcomes}, the term \revise{$\S^2_{iy}(1)  / p_i  + \S^2_{iy}(0)  / (1 - p_i ) - \S_{ i \tau }^2$ is the variance of   $\surd{n_i} (\hat \tau_{i \cdot } - \tau_{i\cdot})$, whose weighted average tends to a finite and positive limit.} It implies that the variance $\Sunadj^2$ is of order $1/N$ because {\color{black}
\begin{equation}
\Sunadj^2 = \sumi c_i \frac{\ni}{N} \Big\{  \frac{\S^2_{iy}(1)}{ \nti } +   \frac{\S^2_{iy}(0)}{ \nci } - \frac{ \S_{ i \tau}^2 }{ \ni }\Big\} = \frac{1}{N} \sumi c_i \Big\{  \frac{\S^2_{iy}(1)}{ p_i} +   \frac{\S^2_{iy}(0)}{ 1 - p_i} - \S_{ i \tau}^2 \Big\}. \nonumber
\end{equation}}
\end{remark}

\begin{theorem}
\label{thm:clt}
If Conditions \ref{cond:prop}--\ref{cond:moments_outcomes} hold, then $( \tauunadj - \tau ) / \Sunadj $  converges in distribution to $\mathcal{N}(0,1)$. Furthermore, if for each stratum $i$,  $2\leq \nti \leq \ni - 2$, then, the  estimator $N \sunadj^2$ converges in probability to the limit of \revise{$  \sumi c_i \big\{ \S^2_{iy}(1)  / p_i  + \S^2_{iy}(0)  / (1 - p_i )  \big\} $}, which is no less than that of $N \Sunadj^2$, and the difference is the limit of $  \sumi c_i \S^2_{i \tau} $.
\end{theorem}

\begin{remark}
Theorem~\ref{thm:clt} is very general in two aspects. \revise{First, the asymptotic normality holds for arbitrary number of strata $B$ and their sizes $n_i$, such that $N \rightarrow \infty$, including the asymptotic regime where $B\rightarrow \infty$ as $N \rightarrow \infty$ with $n_i$ fixed. It includes the paired and finely stratified randomized experiments as special cases. Second, $p_i$ could be different across strata.} A caveat to our analysis, however, is that when estimating the asymptotic variance, we require each stratum to have at least two treated and two control units, which rules out the paired and finely stratified randomized experiments. This is because \revise{$\s^2_{iy}(1) $ or $\s^2_{iy}(0) $} is not well-defined when there is exactly one treated or one control unit. Suitable variance estimators can be found in \revise{\citet{Abadie2008}}, \citet{imai2008variance},  and  \citet{fogarty2018finely, fogarty2018}.
\end{remark}

Theorem~\ref{thm:clt} provides a large-sample Normal approximation of $\tauunadj$, which can produce an asymptotically conservative confidence interval for the average treatment effect \revise{$\tau$:  $[  \tauunadj - q_{\alpha/2} \hat \sigma_{\unadj} ,    \tauunadj + q_{\alpha/2} \hat \sigma_{\unadj} ]$}, where $\alpha$ is the significance level and $q_{\alpha/2}$ is the upper $\alpha/2$ quantile of a standard normal distribution. For comparison, in the special cases of  finely stratified and paired randomized experiments, \citet{fogarty2018finely,fogarty2018} have established the asymptotic normality of $\tauunadj$, but under stronger conditions, the fourth moment conditions on the potential outcomes. Our theorem generalizes their results.

\section{Covariates adjustment}
\subsection{\revise{Covariates adjustment for many small strata}}
Covariates sometimes are predictive for the potential outcomes and we may improve the treatment effect estimation efficiency  by adjusting their  imbalances between treatment and control groups. \revise{Let $B_{ij}$ be the indicator of whether or not unit $j$ is in stratum $i$. Following \citet{lin2013}, we include the treatment-by-covariate interactions in the regression. Precisely, we first compute the weighted potential outcomes and covariates, $y_{ij}^w = y_{ij} w_{ij}$, $\X_{ij}^w = ( \X_{ij}  - \X_{i \cdot} ) w_{ij}$, $B_{ij}^w = (B_{ij} - c_i ) w_{ij}$, where $w_{ij} = Z_{ij} \{ n_i / (\nti - 1) \}^{1/2} + (1 - Z_{ij} ) \{ n_i / ( \nci - 1 ) \}^{1/2}$. Then, we perform the following linear regression:}
\begin{eqnarray}
\label{eqn:weighted_reg}
y_{ij}^w & \sim &  \alpha w_{ij} + \tau Z_{ij} w_{ij} + \sum_{k=2}^{B} B_{kj}^w \gamma_{0k} +  Z_{ij} \sum_{k=2}^{B} B_{kj}^w \gamma_{1k} + \big(   \X_{ij}^w  \big)^\T  \xi_0 +   Z_{ij} \big(   \X_{ij}^w   \big)^\T   \xi_1    .
\end{eqnarray}
The regression-adjusted average treatment effect estimator $\tauols$ is the ordinary least squares estimator of $\tau$. The weights $w_{ij}$ are different from those in \citet{GerberGreen2012}, but they can guarantee the variance reduction. Suppose that $p_i$ is asymptotically the same across strata, then, if one of the following conditions holds, (i) $B$ is bounded and $\ni \rightarrow \infty$, and (ii) $B \rightarrow \infty$ and $n_i$ is bounded and asymptotically the same across strata, the weighted regression \eqref{eqn:weighted_reg} is asymptotically equivalent to the unweighted regression, otherwise, the latter may lose efficiency, while the former does not.  In addition, the former sometimes has smaller variance than the latter when $p_i$'s are different across strata. \revise{Regression \eqref{eqn:weighted_reg} is equivalent to performing two weighted regressions. That is, we compute the weighted least squares estimator based on the treated units,
\begin{equation}
\label{eqn:hat_beta_1}
\hat \beta_1 = \argmin\limits_{\beta}  \sumi \frac{ c_i }{n_{1i} - 1} \sumj Z_{ij}  \Big[  y_{ij}(1) - \hat y_{i\cdot}(1) - \big\{ \X_{ij} - \hat \X_{i\cdot}(1) \big\}^\T \beta \Big]^2. \nonumber
\end{equation}
Similarly, we compute $\hat \beta_0$ based on the control units.} \revise{Denote
\begin{equation}
\label{eqn:hat_tau_i}
\tauiols = \Big[  \hat y_{i\cdot}(1) - \big\{ \hat \X_{i\cdot}(1)  - \X_{i\cdot} \big\}^\T \hat \beta_{1} \Big]  - \Big [ \hat y_{i\cdot}(0) - \big\{ \hat \X_{i\cdot}(0) - \X_{i\cdot} \big\}^\T \hat \beta_{0} \Big], \nonumber
\end{equation}
where $ \{ \hat \X_{i\cdot}(1)  - \X_{i\cdot} \}^\T \hat \beta_{1}  $ and $ \{ \hat \X_{i\cdot}(0) - \X_{i\cdot} \}^\T \hat \beta_{0}$} adjust imbalances of covariate means between treatment and control groups in stratum $i$. The regression-adjusted average treatment effect estimator is the weighted average: \revise{
\[ \tauols    = \sumi c_i \tauiols = \Big[ \hat y_{\cdot \cdot}(1) - \big\{ \hat \X_{\cdot \cdot}(1) - \X_{\cdot \cdot}  \big\}^\T \hat \beta_1  \Big] -  \Big[ \hat y_{\cdot \cdot}(0) - \big\{ \hat \X_{\cdot \cdot}(0) - \X_{\cdot \cdot}  \big\}^\T \hat \beta_0  \Big].
\]}


To study the asymptotic behavior of $\tauols$, we project the potential outcomes onto the space spanned by the linear combination of the covariates with stratum-specific intercepts. Specifically, we  define $\beta_1 $ as the population projection coefficient vector obtained from weighted least squares,
\begin{equation}
\label{eqn:beta_1}
\beta_1 = \argmin\limits_{\beta} \sumi \frac{c_i}{\ni - 1} \sumj \Big[ y_{ij}(1) - y_{i\cdot}(1) - \big\{ \X_{ij} - \X_{i\cdot}  \big\}^\T \beta \Big]^2. \nonumber
\end{equation}
And $\beta_0$ is similarly defined. Then, we define projection errors $\v_{ij}(1)$ and $\v_{ij}(0)$, such that
\begin{equation}
\label{eqn:decomp}
y_{ij}(1) = y_{i\cdot}(1) + \left( \X_{ij} - \X_{i\cdot} \right)^\T \beta_{1} + \v_{ij}(1), \quad y_{ij}(0) = y_{i\cdot}(0) + \left( \X_{ij} - \X_{i\cdot} \right)^\T \beta_{0} + \v_{ij}(0). \nonumber
\end{equation}
All terms in the above equations are fixed population quantities. They are different from linear regression models where $\v_{ij}(1)$ and $\v_{ij}(0)$ are  independently and identically distributed random errors. Both $\v_{ij}(1)$ and $\v_{ij}(0)$ are not observed, but we can estimate them using residuals. The population variance $\S^2_{i \v}(1)$ can be estimated by
\begin{equation}
\hat \s^2_{i \v}(1)  = \frac{1}{\nti - 1} \sumj Z_{ij} \Big[ y_{ij}(1) - \hat y_{i\cdot}(1) - \big\{ \X_{ij} -  \hat \X_{i\cdot}(1) \big\}^\T \hat \beta_{1} \Big]^2.
\end{equation}
Similarly, we define $\hat \s^2_{i\v}(0)$ to estimate  $\S^2_{i \v}(0)$. The variance $  \S^2_{i \v}$ is generally not estimable.  Denote the weighted population variance and its estimator as
\[ \sigma^2_{\ols} = \sumi c_i^2  \Big\{ \frac{ \S^2_{i \v}(1) }{ \nti }  +  \frac{ \S^2_{i \v}(0) }{ \nci } -  \frac{ \S^2_{ i \v} }{ \ni } \Big\}, \quad \hat \sigma^2_{\ols} = \sumi c_i^2  \Big \{ \frac{ \hat \s^2_{i \v}(1) }{ \nti }  + \frac{ \hat \s^2_{i \v}(0) }{ \nci }  \Big \},
\]
respectively. We provide conditions under which covariates adjustment improves or at least does not hurt the precision of average treatment effect estimation in stratified randomized experiments.
\begin{condition}
\label{cond:max_dist_x}
For each covariate $k = 1,\ldots, K$,
$$\frac{1}{N} \maxi \maxj ( \x_{ijk} - \x_{i\cdot k} )^2   \rightarrow 0.$$
\end{condition}
\begin{condition}
\label{cond:variance_e}
The weighted variances $\sumi c_i  \S^2_{i \v}(1)  / p_i $, $ \sumi c_i  \S^2_{i \v}(0)  / (1 - p_i ) $ and $\sumi c_i   \S_{i \v}^2 $ tend to finite limits, positive for the first two, and the limit of $ \sumi c_i \big\{ \S^2_{i \v}(1)  / p_i  +   \S^2_{i \v}(0)  / (1 - p_i )  -   \S_{i \v}^2 \big\} $ is positive.
\end{condition}
\begin{condition}
\label{cond:moments_x}
The weighted covariance matrix $S_{\X \X} = \sumi c_i  \S_{i \X \X}$ converges to a finite, invertible matrix; the weighted covariances, $S_{\X y}(1) = \sumi c_i \S_{i \X y}(1)$ and $S_{\X y}(0) = \sumi c_i \S_{i \X y}(0) $, and the weighted absolute covariances,  $\sumi c_i |  \S_{i \X \X} |$,  $ \sumi c_i | \S_{i \X y}(1) | $, and $\sumi c_i |
\S_{i \X y}(0) |$, converge to finite limits.
\end{condition}

\begin{theorem}
\label{thm:clt_ols}\revise{
If Conditions \ref{cond:prop}--\ref{cond:moments_x} hold and assume that for each stratum $i$,  $2\leq \nti \leq \ni - 2$, then, $(\tauols - \tau) / \sigma_{\ols}$ converges in distribution to $\mathcal{N}(0,1)$ as $N \rightarrow \infty$.  Furthermore, if $p_{i,\infty} = p $ for $i=1,\ldots,B$, the difference between the asymptotic variance of $\surd{N} \tauols$ and $\surd{N} \tauunadj$ is the limit of $ - \sumi c_i \Delta_i^2$, where
$$
\Delta_i^2 =\big( \beta^i \big)^\T \S_{i\X \X} \big( \beta^i \big)  /   \big\{ p_i  ( 1 - p_i ) \big\}, \quad  \beta^i  = (1 - p_i ) \beta_{1} + p_i \beta_{0}.
$$
The estimator $ N \hat \sigma^2_{\ols}  $ converges in probability to the limit of $  \sumi c_i \{ \S^2_{i \v}(1)  / p_i + \S^2_{i \v}(0)  / ( 1 - p_i )  \} $, which is no less than that of $N \sigma^2_{\ols}$ and no greater than that of  $N \sunadj^2$. The differences, $ N (  \hat \sigma^2_{\ols}  -  \sigma^2_{\ols} ) $ and $ N ( \hat \sigma^2_{\ols} -  \sunadj^2 ) $, converge in probability to  the limits of $  \sumi c_i \S^2_{i \v } $ and $ -   \beta_1^\T  \S_{\X\X} \beta_1/p - \beta_0^\T \S_{\X\X} \beta_0/(1-p) $, respectively. }
\end{theorem}

\begin{remark}
As pointed out by \citet{cochran1977}, stratified random sampling does not always result in a smaller variance for the population mean than that given by simple random sampling. However, if intelligently used, for example, the sampling proportion is the same in all strata, stratified random sampling gives a variance no larger than that of the simple random sampling. For the asymptotic normality, \revise{both $B$ and $\ni$ can tend to infinity, and $p_i$ can be different across strata. Our main requirement for improving efficiency is that $p_i$ tends uniformly to a common limit $p$ $(0 < p < 1)$ as $N \rightarrow \infty$. In order to define $\hat \beta_1$ and $\hat \beta_0$, and to obtain an appropriate variance estimator, we assume that each stratum has at least two treated and two control units.}
\end{remark}

\begin{remark}\revise{
If we use fixed adjusted-vectors, such as $\beta_1$ and $\beta_0$, in the definition of $\tauols$, it is unbiased. However, since the adjusted-vectors $\hat \beta_1$ and $\hat \beta_0$ are obtained from the data,  $\tauols$ introduces a small-sample bias, but the bias diminishes rapidly. Under stratum-specific fourth moment conditions on the covariates and potential outcomes, both $\hat \beta_1 - \beta_1$ and $\hat \X_{\cdot \cdot}(1) - \X_{\cdot \cdot}$ are $\surd{N}$-consistent. The bias of $\tauols$ is of order $1/N$.}
\end{remark}

Theorem \ref{thm:clt_ols} implies that the regression-adjusted average treatment effect estimator $\tauols$ is consistent and asymptotically normal with asymptotic variance no more than that of the stratified difference-in-means estimator $\tauunadj$. Thus, covariates adjustment improves or does not hurt precision, at least asymptotically. The improvement depends on whether  the covariates are predictive to the potential outcomes in such a way that $\Delta_i^2 > 0$ for some stratum $i$. This theorem also provides a conservative variance estimator, which is consistent  if and only if the stratum-specific treatment effect is constant asymptotically, or equivalently, the limit of $\sumi c_i \S^2_{i \v } $ is zero. Moreover, the variance estimator $\hat \sigma^2_{\ols}$ is asymptotically no worse than the unadjusted variance estimator $ \sunadj^2$, and thus, we can use it to construct a large-sample confidence interval for the average treatment effect that is at least as good as the one constructed by $  \sunadj^2$.

\subsection{Covariates adjustment for a few large strata}

In some \revise{stratified randomized experiments}, the number of strata $B$ is small and each has a large  stratum size $n_i$. For example, in a clinical trial with many males and females, experimental units are stratified according to their gender, resulting in two large strata.  The stratified randomized experiment is hence the independent combination of a few completely randomized experiments, each with a large sample size. It should be more efficient to use stratum-specific adjusted coefficients $\hat \beta_{1i}$ and $\hat \beta_{0i}$ instead of common ones for all strata. More precisely, within stratum $i$, we regress the outcomes $y_{ij} \ (j=1,\ldots,\ni)$, on covariates in the treated and control \revise{groups} separately, and obtain  $\hat \beta_{1i}$ from the ordinary least squares, \revise{
\begin{equation}
\hat \beta_{1i}  = \argmin\limits_{\beta} \sumj Z_{ij} \Big [ y_{ij}(1) - \hat y_{ i \cdot}(1) -  \big\{ X_{ij} - \hat X_{i\cdot}(1) \big\}^\T \beta \Big ]^2. \nonumber
\end{equation}
And $\hat \beta_{0i}$ is similarly defined.} The above procedure is equivalent to including the interaction of stratum indicator and covariates in the regression, hopefully resulting in smaller variance of approximation errors, and thus, decreasing the variance of average treatment effect estimator. The stratum-specific average treatment effect estimator is \revise{
\begin{equation}
\tauiolsint = \Big[  \hat y_{i\cdot}(1) - \big\{ \hat X_{i\cdot}(1)  - X_{i\cdot} \big\}^\T \hat \beta_{1i} \Big]  - \Big[ \hat y_{i\cdot}(0) - \big\{ \hat y_{ i \cdot}(0) - X_{i\cdot} \big\}^\T \hat \beta_{0i} \Big]. \nonumber
\end{equation}}
The regression with interaction-adjusted average treatment effect estimator is
\[ \tauolsint = \sumi c_i \tauiolsint . \]

The asymptotic behavior of $\tauolsint$ depends on the projection of the potential outcomes onto the space spanned by the linear combination of covariates, stratum indicator, and their interactions. \revise{We define $\beta_{1i}$  as the projection coefficient vector obtained from ordinary least squares,
\begin{equation}
\beta_{1i} = \argmin\limits_{\beta} \frac{1}{ \ni - 1} \sumj \Big[ y_{ij}(1) - y_{i\cdot}(1) - \big\{ X_{ij} - X_{i\cdot} \big\}^\T \beta \Big]^2. \nonumber
\end{equation}
We define $\beta_{0i}$ similarly. Then, we define projection errors $\e_{ij}(1)$ and $\e_{ij}(0)$, such that
\begin{equation}
\label{eqn:decomp_inter}
y_{ij}(1) = y_{i\cdot}(1) + \left( X_{ij} - X_{i\cdot} \right)^\T \beta_{1i} +  \e_{ij}(1), \quad y_{ij}(0) = y_{i\cdot}(0) + \left( X_{ij} - X_{i\cdot} \right)^\T \beta_{0i} +  \e_{ij}(0).   \nonumber
\end{equation}
Again, all quantities in the above projections are not random. To estimate the asymptotic variance of $\tauolsint$, we define
\[
\hat \s^2_{i \e}(1)  = \frac{1}{\ni - K - 1} \sumj  Z_{ij} \Big[ y_{ij}(1) - \hat y_{i\cdot}(1) - \big\{ X_{ij} -  \hat X_{i\cdot}(1) \big\}^\T \hat \beta_{1i} \Big]^2,
\]
and we define $\hat \s^2_{i \e}(0)$ similarly. Following \citet{bloniarz2015lasso}, we adjust  the degrees of freedom. The variance of $\tauolsint$ and its estimator are
\[
\sigma^2_{\olsint} = \sumi c_i^2  \Big\{  \frac{ \S^2_{i\e}(1) }{ \nti }  +  \frac{ \S^2_{i\e}(0) }{ \nci }  -  \frac{ \S^2_{i \e} }{ \ni } \Big\}, \quad \hat \sigma^2_{\olsint} = \sumi c_i^2  \Big\{ \frac{ \hat  \s^2_{i\e}(1) }{ \nti }  +  \frac{ \hat  \s^2_{i\e}(0) }{
\nci }  \Big \},
\]
respectively.} We now present the condition and asymptotic result as follows.

\begin{condition}
\label{cond:large_ni}
(a) There exists a constant $B_{\max}$ not depending on $N$, such that $B \leq B_{\max}$ for every $N$. The stratum size $n_i \rightarrow \infty$ as $N \rightarrow \infty$ for $i = 1,\ldots, B$. (b) The stratum-specific covariance matrix $\S_{i\X\X}$ converges to a finite, invertible matrix. The stratum-specific variances $\S^2_{iy}(1) $, $ \S^2_{iy}(0) $, $\S^2_{i\tau } $, $\S^2_{i\e}(1) $, $ \S^2_{i\e}(0) $, $\S^2_{i \e}$, and the covariances $\S_{i\X y}(1)$, $\S_{i\X y}(0)$ converge to finite limits.
\end{condition}

\begin{theorem}
\label{thm:clt_ols_int}\revise{
If Conditions \ref{cond:prop}, \ref{cond:max_dist_y}, \ref{cond:max_dist_x} and \ref{cond:large_ni} hold, then, $(\tauolsint - \tau) / \sigma_{\olsint}$ converges in distribution to $\mathcal{N}(0,1)$ as $N \rightarrow \infty$. Furthermore, the difference between the asymptotic variance of $\surd{N} \tauolsint$ and $ \surd{N} \tauunadj$ tends to the limit of $ - \sumi c_i \tilde \Delta_i^2 \leq 0 $, where 
$$ \tilde \Delta_i^2 =  \big( \tilde \beta^i \big)^\T \S_{i\X \X} \big( \tilde \beta^i \big) / \big\{ p_i  ( 1 - p_i ) \big\}, \quad \tilde \beta^i  = (1 - p_i ) \beta_{1i} + p_i \beta_{0i}. $$
The estimator  $ N \hat \sigma^2_{\olsint}  $ converges in probability to the limit of $  \sumi c_i \{ \S_{i \e}^2(1)  / p_i + \S_{i \e}^2(0)  / ( 1 - p_i )  \} $, which is no less than that of $N \sigma^2_{\olsint}$ and no greater than that of  $N \sunadj^2$. The differences, $ N (  \hat \sigma^2_{\olsint}  -  \sigma^2_{\olsint} ) $ and $N (  \hat \sigma^2_{\olsint} -  \sunadj^2 ) $, converge in probability to the limits of $  \sumi c_i \S^2_{i \e } $ and $ -    \sumi (c_i /p_i ) \beta_{1i}^\T \S_{i\X\X} \beta_{1i} -  \sumi \{ c_i /( 1 - p_i ) \} \beta_{0i}^\T \S_{i\X\X} \beta_{0i} $, respectively.}
\end{theorem}

\begin{remark}
\revise{The asymptotic normality result is a direct extension of \citet{fcltxlpd2016}}. When the \revise{projections are} homogeneous across strata, $\beta_{1i} = \beta_1$ and $\beta_{0i} = \beta_0$, for $i= 1,\ldots, B$,  $\tauolsint$ and $\tauols$ are asymptotically equivalent. However, when the \revise{projections are} heterogeneous, $\tauolsint$ is more efficient  in general than $\tauols$; see the following Corollary~\ref{col:var_ols_int}.
\end{remark}

\begin{remark}
\revise{
When $n_i \rightarrow \infty$, the bias of $\hat \tau_{i, \olsint}$ is of order $1/ n_i$ \citep{lin2013}. Therefore, when $B$ is bounded and $n_i \rightarrow \infty$, the bias of $\tauolsint$ is of order $1/N$, which tends to zero much faster than the variance.}
\end{remark}

Theorem \ref{thm:clt_ols_int} states that, compared with the stratified difference-in-means estimator, $\tauunadj$,  the regression with interaction-adjusted  estimator  $\tauolsint$ improves or  does not hurt precision when there are only a few large strata. The asymptotic variance can be conservatively estimated by the weighted residual sum of squares, \revise{$\hat \sigma^2_{\olsint}$}, which performs no worse than the unadjusted estimator \revise{$ \sunadj^2$}. Moreover, this theorem does not require asymptotically the same proportion of treated units across strata as assumed in Therorem~\ref{thm:clt_ols}.

\begin{corollary}
\label{col:var_ols_int}
Under Conditions~\ref{cond:prop}--\ref{cond:large_ni},  if for each stratum $i$, $2\leq \nti \leq \ni - 2$ and \revise{$p_{i,\infty} = p $}, then the following holds, in probability for the second one,
\[
\revise{N \sigma^2_{\olsint} \leq N \sigma^2_{\ols} \leq N \Sunadj^2, \quad N \hat \sigma^2_{\olsint} \leq N  \hat \sigma^2_{\ols} \leq N \sunadj^2.}
\]
\end{corollary}

Corollary~\ref{col:var_ols_int} shows that  the regression with interaction-adjusted estimator $\tauolsint$ performs better (at least no worse) than the regression-adjusted estimator $\tauols$, when there are a few large strata. Furthermore, the asymptotic variance estimator $\hat \sigma^2_{\olsint}$ performs  at least as good as $ \hat \sigma^2_{\ols}$.

\section{Simulation study}
In this section, we conduct a simulation study to evaluate the finite-sample performance of covariates adjustment in stratified randomized experiments. Our data generating process is similar to that of \citet{bloniarz2015lasso}. Specifically, we generate the potential outcomes from the following nonlinear models:
\[ y_{ij}(1) =  \X_{ij}^\T \beta_{11} + \exp{ \left( \X_{ij}^\T \beta_{12} \right) } + \v_{ij}(1) \quad (i = 1,\ldots, B; \ j = 1,\ldots, \ni),
\]
\[ y_{ij}(0) =  \X_{ij}^\T \beta_{01} + \exp{ \left( \X_{ij}^\T \beta_{02} \right) } + \v_{ij}(0) \quad (i = 1,\ldots, B; \ j = 1,\ldots, \ni),
\]
where $\v_{ij}(1)$ and $\v_{ij}(0)$ are generated from Gaussian distribution with mean zero and variance $\sigma^2$ such that the signal-to-noise ratio equals one. The $\X_{ij}$ is a 10-dimensional \revise{covariate}  vector generated from a multivariate normal distribution with mean zero and covariance matrix $\Sigma$: $
\Sigma_{kk} = 1$ and $ \Sigma_{kl} = \rho^{| k - l |}, \ k \neq l, \  k, l = 1,\ldots, 10$, where $\rho$ controls the correlation between covariates and it takes \revise{values} of $0$ and $0.5$. The number of strata, stratum size, and coefficients  $\beta_{11}$, $\beta_{12}$, $\beta_{01}$, $\beta_{02} $ are generated in four scenarios.

(1) Scenario 1. There are many small strata. We set the number of strata $B = 25, 50, 100$, and the stratum size equals ten. The elements of coefficients $\beta_{11}$, $\beta_{12}$, $\beta_{01}$, $\beta_{02} $ are generated as follows: 
$$ \beta_{11k} \sim  t_3, \quad \beta_{12k} \sim 0.1*t_3, \quad  \beta_{01k} \sim \beta_{11k} + t_3, \quad  \beta_{02k} \sim \beta_{12k} + 0.1*t_3, \quad  k=1,\ldots,10, $$
where $t_3$ \revise{denotes} the $t$ distribution with three degrees of freedom.

(2) Scenario 2. There are two large homogeneous strata, $B = 2$, and we consider three values of stratum size, $\ni = 100, 200, 500$. The coefficients  $\beta_{11}$, $\beta_{12}$, $\beta_{01}$, $\beta_{02} $ are generated in the same way as in Scenario 1.

(3) Scenario 3. There are two large heterogeneous strata with same stratum size as in Scenario 2. The difference is that the coefficients  $\beta_{11}$, $\beta_{12}$, $\beta_{01}$, $\beta_{02} $ are generated separately and independently for different stratum.

(4) Scenario 4. There are two large heterogeneous strata of size $100$, the same as Scenario 3, and $B=25, 50, 100$ small strata of size $10$, the same as Scenario 1. The covariate vector, coefficients, and potential outcomes are generated in the same way as the first three scenarios, except that the number of covariates is three.

Both the covariates and potential outcomes are generated once, and then kept fixed. Thereafter, we simulate a stratified randomized experiment for $10000$ times, assigning $\nti = \ni/2$ units to the treatment and the remainders to the control.  In Scenario 1, since the stratum size is small, we estimate and infer the average treatment effect using only the stratified difference-in-means estimator $\tauunadj$ and the regression-adjusted estimator $\tauols$, while in Scenarios 2 and 3, we add the regression with interaction-adjusted estimator $\tauolsint$. We compare their \revise{performances} in terms of absolute bias, standard deviation, root mean squared error, empirical coverage probability and mean confidence interval length of $95\%$ confidence interval constructed as \revise{$[ \hat \tau - 1.96  \hat \sigma, \hat \tau + 1.96 \hat \sigma ]$}, where $\hat \tau$ is the average treatment effect estimator and \revise{$\hat \sigma$} is the estimated standard deviation.

\begin{figure}
\includegraphics[width=\textwidth,height=.8\textheight]{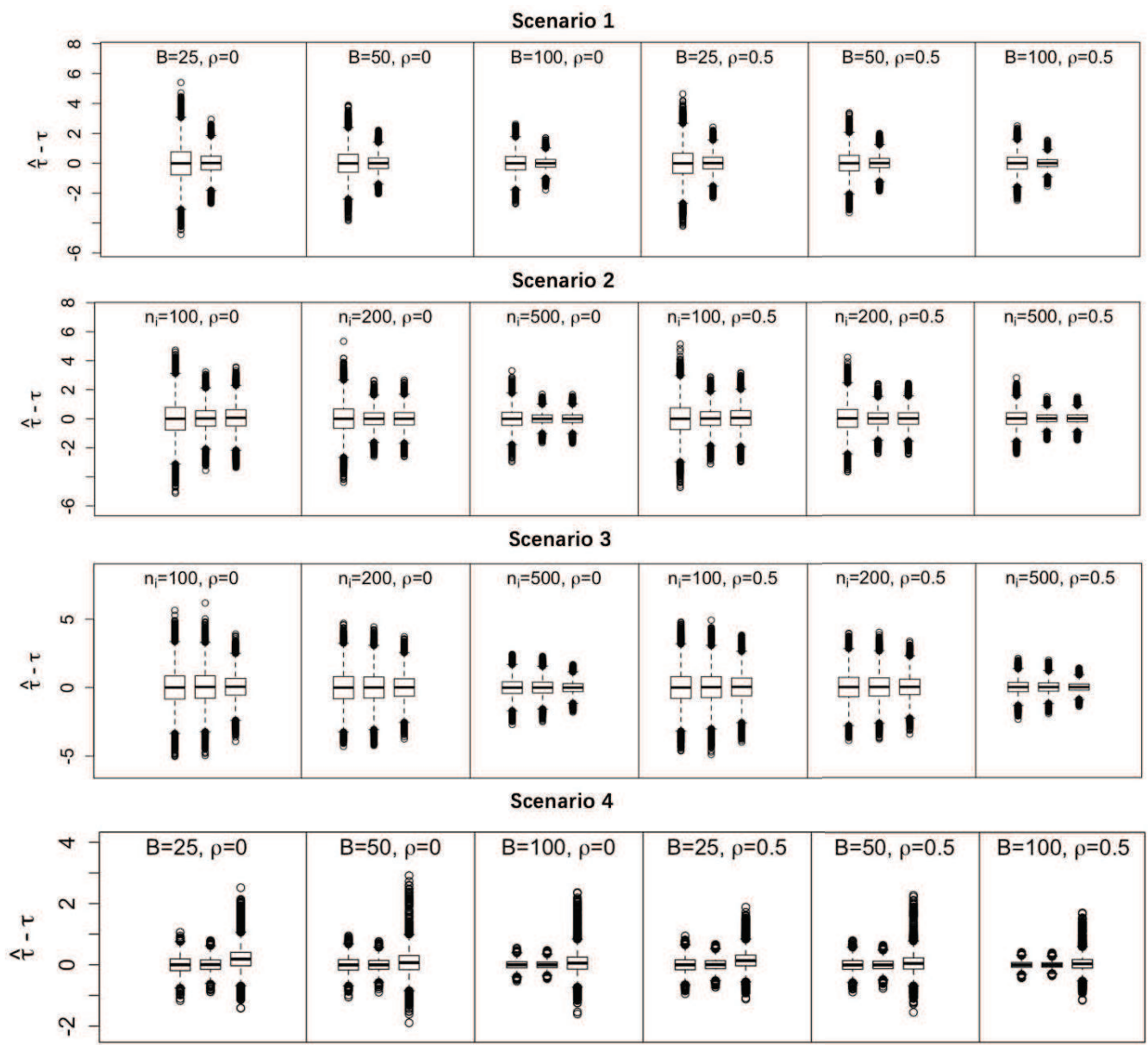}
\caption{Box plot of average treatment effect estimator minus the true average treatment effect, $\hat \tau - \tau$. In each sub-slot of each scenario, the box \revise{plots correspond to} the methods, from left to right, $\tauunadj$,  $\tauols$, and $\tauolsint$. In Scenario 1, we do not compute and present $\tauolsint$ \revise{because both the number of treated and control units are smaller than the number of covariates}.}
\label{fig:box}
\end{figure}

The results are shown in Figure~\ref{fig:box} and Tables~\ref{tab:tab1}--\ref{tab:tab3}.  In all scenarios, the bias of each method is substantially smaller, more than 10 times, than its standard deviation. This is because $\tauunadj$ is unbiased and both $\tauols$ and $\tauolsint$ are asymptotically unbiased.  In the first two scenarios, $\tauols$ performs the best, reducing the standard deviation and root mean squared error of $\tauunadj$ by  $32\%  - 43\%$, although it introduces a small and negligible amount of finite-sample bias. The estimator $\tauols$ performs slightly better than $\tauolsint$ \revise{in Scenario 2 where the potential outcomes are generated homogeneously across  strata}. In Scenario 3 where there are two large heterogeneous  strata, $\tauols$ can still reduce the standard deviation and root mean squared error of $\tauunadj$ by $3\% - 10\%$, although it is not as significant as in the first two scenarios. In this scenario, $\tauolsint$ performs the best, which further reduces the standard deviation and root mean squared error of $\tauols$ by $17 \% -  33\%$.  In Scenario 4, the regression-adjusted estimator $\tauols$ performs the best. It decreases the root mean squared error and mean confidence interval length by $9\% - 20\%$ and $5\% - 8\%$, respectively, when compared with the unadjusted estimator $\tauunadj$. The regression with interaction-adjusted estimator $\tauolsint$ does not work at all, which is even not consistent. This is because the number of observations in small strata is too small to accurately estimate the regression coefficient. Therefore, $\tauolsint$ should be used with caution when there are many small strata. Furthermore, for all scenarios and methods, the variance estimators are conservative, resulting in conservative confidence intervals with empirical coverage probabilities much higher than the confidence level. However, when compared with the unadjusted estimator $\tauunadj$, $\tauols$ substantially reduces the mean confidence interval length, $4\%  - 19\%$, in all scenarios. The $\tauolsint$ performs similarly to $\tauols$ in Scenario 2 and performs better, $5 \%  - 10 \%$, than $\tauols$ in Scenario 3.

\begin{table}
\centering
 \caption{\label{tab:tab1} Comparison results of average treatment effect estimates in Scenario 1}
\begin{tabular*}{\hsize}{@{\extracolsep{\fill}}cccccccc}
  \hline
$\rho$ & $B$ & methods & Bias $(\times 1000)$ & SD $(\times 100)$ & $\surd \textnormal{MSE}$ $(\times 100)$ & CP $(\%)$ & CI length $(\times 100)$ \\ 
  \hline
0 & 25 & $\tauunadj$ & 2(6) & 114(7) & 114(23) & 99.0 & 232 \\ 
& & $\tauols$  & 14(2) & 68(5) & 68(8) & 99.6 & 190 \\  \cline{2-8}
& 50 & $\tauunadj$  & 2(3) & 89(6) & 89(14) & 98.8 & 200 \\ 
& & $\tauols$ & 4(1) & 52(3) & 52(5) & 99.5 & 164 \\  \cline{2-8}
& 100 & $\tauunadj$  & 0(2) & 66(4) & 66(8) & 98.8 & 173 \\ 
& & $\tauols$ & 1(1) & 38(3) & 38(3) & 99.7 & 142 \\  \hline
0.5 & 25& $\tauunadj$  & 1(4) & 99(7) & 99(18) & 98.9 & 216 \\ 
& & $\tauols$ & 13(1) & 57(4) & 57(6) & 99.6 & 174 \\  \cline{2-8}
& 50 & $\tauunadj$  & 1(3) & 77(5) & 77(10) & 98.8 & 187 \\ 
& & $\tauols$ & 3(1) & 46(3) & 46(4) & 99.5 & 154 \\  \cline{2-8}
& 100 & $\tauunadj$  & 1(2) & 59(4) & 59(6) & 98.8 & 164 \\ 
& & $\tauols$ & 1(1) & 34(2) & 34(2) & 99.6 & 135 \\ 
   \hline
\end{tabular*}
SD, standard deviation; $\surd \textnormal{MSE}$, root mean squared error; CP, coverage probability; CI length, mean confidence interval length. The numbers in parentheses are the corresponding standard errors estimated by using the bootstrap with $500$ replications.
\end{table}

\begin{table}
\centering
 \caption{\label{tab:tab2} Comparison results of average treatment effect estimates in Scenario 2 and 3}
\begin{tabular*}{\hsize}{@{\extracolsep{\fill}}cccccccc}
  \hline
  $\rho$ & $\ni$ & methods & Bias $(\times 1000)$ & SD $(\times 100)$ & $\surd \textnormal{MSE}$ $(\times 100)$ & CP $(\%)$ & CI length $(\times 100)$ \\ 
  \hline
  &&&& Scenario 2 &&& \\
  \hline
  0 & 100 & $\tauunadj$  & 4(6) & 116(8) & 116(23) & 99.1 & 234 \\ 
 & & $\tauols$ & 26(3) & 78(5) & 78(11) & 99.3 & 196 \\ 
 & & $\tauolsint$ & 67(3) & 83(5) & 83(12) & 99.2 & 203 \\  \cline{2-8}
 & 200 & $\tauunadj$  & 4(4) & 99(7) & 99(17) & 98.8 & 213 \\ 
 & & $\tauols$  & 4(2) & 61(4) & 61(6) & 99.5 & 177 \\ 
 & & $\tauolsint$ & 1(2) & 63(4) & 63(7) & 99.4 & 180 \\  \cline{2-8}
 & 500 & $\tauunadj$ & 1(2) & 66(4) & 66(8) & 98.8 & 173 \\ 
 & & $\tauols$  & 1(1) & 38(3) & 38(3) & 99.6 & 142 \\ 
 & & $\tauolsint$ & 0(1) & 38(2) & 38(3) & 99.6 & 143 \\  \hline
0.5 & 100 & $\tauunadj$  & 5(6) & 110(7) & 110(21) & 98.9 & 226 \\ 
 & & $\tauols$  & 22(2) & 70(5) & 70(8) & 99.3 & 185 \\ 
 & & $\tauolsint$ & 59(3) & 74(5) & 74(9) & 99.2 & 191 \\  \cline{2-8}
   & 200 & $\tauunadj$  & 5(4) & 91(6) & 91(14) & 98.8 & 204 \\ 
 & & $\tauols$  & 3(1) & 56(4) & 56(5) & 99.5 & 170 \\ 
 & & $\tauolsint$ & 1(2) & 58(4) & 58(6) & 99.4 & 173 \\  \cline{2-8}
 & 500 & $\tauunadj$ & 0(2) & 59(4) & 59(6) & 98.9 & 164 \\ 
 & & $\tauols$  & 1(1) & 34(2) & 34(2) & 99.6 & 135 \\ 
 & & $\tauolsint$ & 0(1) & 34(2) & 34(2) & 99.6 & 135 \\  
  \hline
  &&&& Scenario 3 &&& \\
  \hline
0 & 100 & $\tauunadj$  & 8(7) & 124(8) & 124(28) & 99.1 & 244 \\ 
 & & $\tauols$  & 44(6) & 121(8) & 121(26) & 98.6 & 231 \\ 
 & & $\tauolsint$ & 58(4) & 91(6) & 91(15) & 99.3 & 213 \\  \cline{2-8}
 & 200 & $\tauunadj$  & 1(6) & 117(8) & 117(22) & 98.9 & 231 \\ 
 & & $\tauols$  & 7(5) & 112(8) & 112(20) & 98.6 & 221 \\ 
 & & $\tauolsint$ & 7(4) & 92(6) & 92(14) & 99.1 & 208 \\ \cline{2-8}
 & 500 & $\tauunadj$  & 1(2) & 63(4) & 63(7) & 99.1 & 172 \\ 
 & & $\tauols$  & 0(2) & 58(4) & 58(6) & 99.1 & 166 \\ 
 & & $\tauolsint$ & 1(1) & 43(3) & 43(3) & 99.5 & 150 \\  \hline
0.5 & 100 & $\tauunadj$  & 7(6) & 117(8) & 117(24) & 99.3 & 240 \\ 
 & & $\tauols$  & 34(6) & 113(8) & 113(23) & 98.8 & 227 \\ 
 & & $\tauolsint$ & 48(4) & 97(6) & 97(16) & 99.1 & 215 \\  \cline{2-8}
 & 200 & $\tauunadj$  & 1(4) & 102(7) & 102(17) & 99.0 & 217 \\ 
 & & $\tauols$  & 7(4) & 96(6) & 96(15) & 98.8 & 207 \\ 
 & & $\tauolsint$ & 7(3) & 83(5) & 83(11) & 99.1 & 197 \\  \cline{2-8}
 & 500 & $\tauunadj$  & 0(1) & 50(3) & 50(4) & 99.2 & 157 \\ 
 & & $\tauols$  & 1(1) & 45(3) & 45(3) & 99.3 & 149 \\ 
 & & $\tauolsint$ & 2(0) & 33(2) & 33(2) & 99.6 & 134 \\ 
  \hline
\end{tabular*}
SD, standard deviation; $\surd \textnormal{MSE}$, root mean squared error; CP, coverage probability; CI length, mean confidence interval length. The numbers in parentheses are the corresponding standard errors estimated by using the bootstrap with $500$ replications.
\end{table}

\begin{table}
\centering
 \caption{\label{tab:tab3} Comparison results of average treatment effect estimates in Scenario 4}
\begin{tabular*}{\hsize}{@{\extracolsep{\fill}}cccccccc}
  \hline
$\rho$ & $B$ & methods & Bias $(\times 1000)$ & SD $(\times 100)$ & $\surd \textnormal{MSE}$ $(\times 100)$ & CP $(\%)$ & CI length $(\times 100)$ \\ 
  \hline
0 & 25 & $\tauunadj$ & 1(1) & 28(3) & 28(2) & 98.9 & 115 \\ 
& & $\tauols$ & 5(1) & 23(3) & 23(2) & 99.3 & 106 \\ 
 & & $\tauolsint$  & 195(3) & 36(5) & 41(9) & 87.5 & 96 \\ \cline{2-8}
 & 50 & $\tauunadj$ & 4(1) & 25(3) & 25(2) & 99.2 & 111 \\ 
& & $\tauols$ & 3(1) & 21(3) & 21(1) & 99.3 & 103 \\ 
 & & $\tauolsint$  & 74(4) & 38(6) & 39(13) & 88.3 & 93 \\ \cline{2-8}
 & 100 & $\tauunadj$ & 0(0) & 14(2) & 14(1) & 99.5 & 85 \\ 
& & $\tauols$ & 1(0) & 13(2) & 13(1) & 99.4 & 80 \\ 
 & & $\tauolsint$  & 56(3) & 33(5) & 34(9) & 78.8 & 73 \\ \cline{2-8}
0.5 & 25 & $\tauunadj$ & 1(1) & 24(3) & 24(2) & 99.0 & 106 \\ 
& & $\tauols$ & 3(1) & 19(2) & 19(1) & 99.3 & 98 \\ 
 & & $\tauolsint$  & 146(2) & 29(4) & 32(6) & 89.9 & 89 \\ \cline{2-8}
 & 50 & $\tauunadj$ & 3(1) & 21(3) & 21(1) & 99.3 & 103 \\ 
& & $\tauols$ & 2(0) & 18(2) & 18(1) & 99.5 & 97 \\ 
 & & $\tauolsint$  & 50(2) & 30(5) & 31(8) & 90.3 & 87 \\ \cline{2-8}
 & 100 & $\tauunadj$ & 1(0) & 12(1) & 12(0) & 99.3 & 76 \\ 
& & $\tauols$ & 1(0) & 11(1) & 11(0) & 99.2 & 72 \\ 
 & & $\tauolsint$  & 38(1) & 24(4) & 24(5) & 81.3 & 64 \\ 
    \hline
\end{tabular*}
SD, standard deviation; $\surd \textnormal{MSE}$, root mean squared error; CP, coverage probability; CI length, mean confidence interval length. The numbers in parentheses are the corresponding standard errors estimated by using the bootstrap with $500$ replications.
\end{table}

Overall, the regression-adjusted estimator $\tauols$ always \revise{improves} the precision of average treatment effect estimation, when compared with the unadjusted estimator $\tauunadj$. Only when there are a few large heterogeneous strata,  regression with interaction-adjusted estimator $\tauolsint$ is preferable. Moreover, the variance estimators are all conservative. These findings agree with our theoretical results.

\section{Analysis of iron deficiency and schooling attainment data}
We \revise{apply covariates adjustment} methods to analyse iron deficiency and schooling attainment data from a stratified randomized trial conducted by \citet{chong2016}. The trial aimed to study the effect of iron deficiency on educational attainment and cognitive ability.

We provide here only a brief description of the experimental setup and refer the reader to \citet{chong2016} for more details. The experiment was carried out in a rural secondary school in the Cajamarca district of Peru between October and December in 2009. The experimental units consist of 215 students who were stratified by the number of years of secondary school completed. There were five strata in total and the strata \revise{sizes range} from 30 to 58. \revise{The researchers conducted a stratified randomized experiment. Within each stratum, approximately one third of the students were randomly assigned to watch an educational video in which a well-known soccer player explained the importance of iron for maximizing energy and encouraged iron supplementation, one third assigned to watch a similar educational video in which the soccer player was replaced by a doctor who encouraged iron supplements for overall health, and one third assigned to watch a `placebo' video unrelated to iron in which a dentist encouraged good oral hygiene.} As suggested by \citet{chong2016}, we \revise{group} the first two videos related to iron as treatment and the placebo video as control. \revise{The proportion of treated units in stratum $i$, $p_i$, is approximately two third, $0.688$, $0.672$, $0.652$, $0.636$, and $0.667$ for the five strata, respectively.} \citet{chong2016} examined the effect of \revise{the} treatment on a variety of outcomes, among them the most important ones are: number of iron supplement pills taken between October and December in 2009; average grade point of the last two quarter; and cognitive ability measured by the average score on five Wii games, {\em Big Brain Academy: Wii Degree}, testing the ability of identification, memorization, analysis, computation, and visualization, respectively. There were a few baseline covariates, and among them we select five related ones as illustrations to perform covariates adjustment. The selected covariates are gender, age, distance to school, first quarter grades, and year of mother's education. 

\begin{figure}
\includegraphics[width=\textwidth,height=.3\textheight]{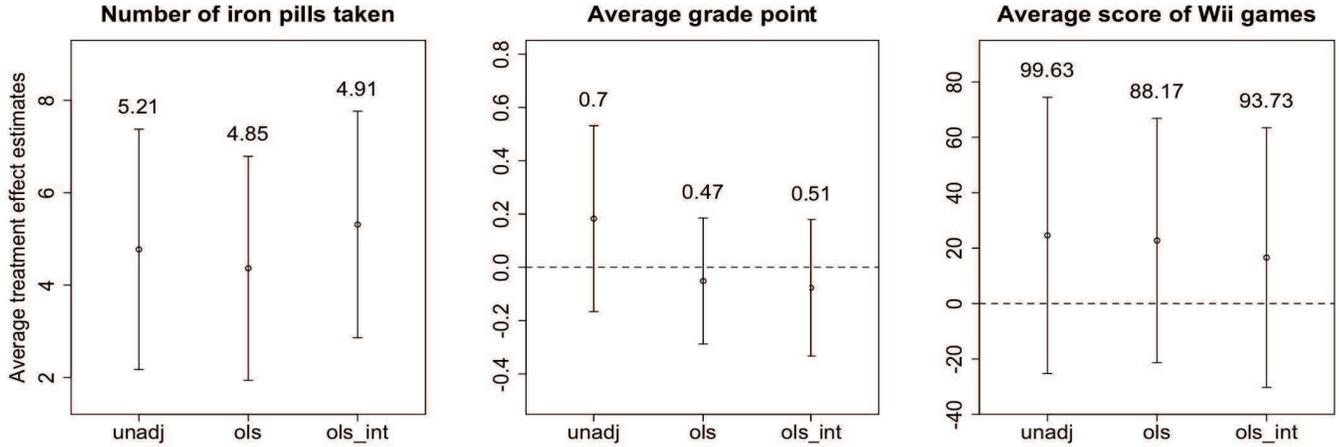}
\caption{Average treatment effect estimates and $95\%$ confidence intervals for three outcomes: number of iron pills taken, average grade score, and average score of Wii games. The circle dots are average treatment effect estimators and the bars are $95\%$ confidence intervals, with the lengths shown on top.}
\label{fig:real}
\end{figure}

\begin{figure}
\includegraphics[width=\textwidth,height=.6\textheight]{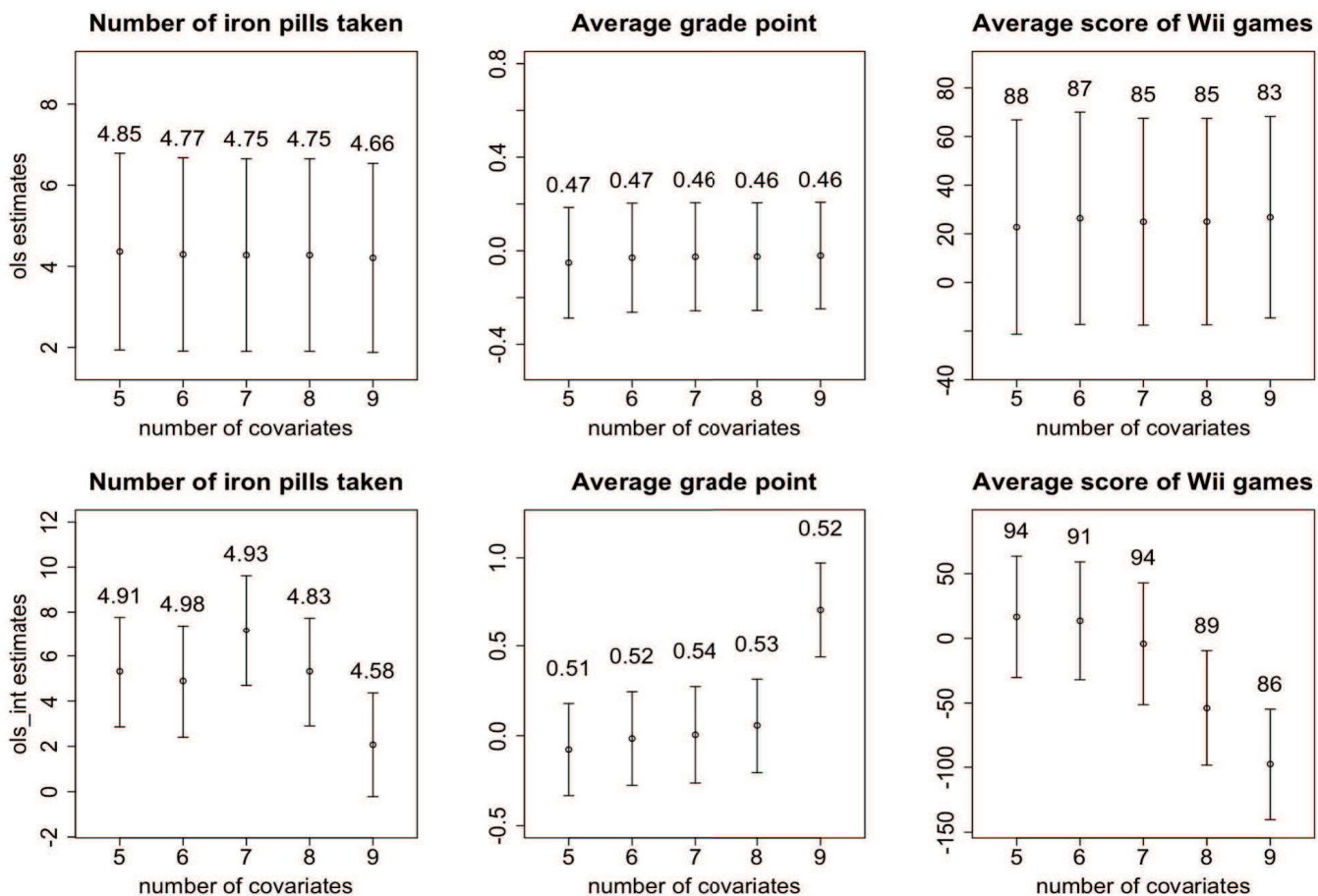}
\caption{Changes of average treatment effect estimates (circle dots) and $95\%$ confidence intervals (bars) for three outcomes: number of iron pills taken, average grade score, and average score of Wii games, when the number of covariates changes from five to nine. The confidence interval lengths are shown on top of the bars. The sub-figures in the first line are the results of $\tauols$, and the sub-figures in the second line are the results of $\tauolsint$.}
\label{fig:real_change_k}
\end{figure}

We estimate and construct $95\%$ confidence intervals for the average treatment effect of the educational video about the importance of iron on the outcomes mentioned above, using the simple stratified difference-in-mean estimator $\tauunadj$,  regression-adjusted estimator $\tauols$, and regression with interaction-adjusted estimator $\tauolsint$. Figure~\ref{fig:real} shows the results. For the first outcome, all three methods produce confidence intervals not containing zero, meaning that the iron educational video has a significant effect on the number of iron supplement pills taken. This conclusion agrees with that of \citet{chong2016}. However, for the other two outcomes related to schooling attainment, the average treatment effects are not significant. \revise{These are not} in  conflict with the findings in \cite{chong2016}, because we infer here the average treatment effect on all experimental units, while \citet{chong2016} studied the effect only on students who are anaemic. Furthermore, compared with $\tauunadj$, both $\tauols$ and $\tauolsint$  produce much shorter confidence intervals, $6 \%  - 32 \%$ shorter, which \revise{are} in accordance with our theoretical results. The $\tauols$ performs the best in this empirical study, whose mean confidence interval lengths are $1\% - 8\%$ shorter than \revise{those} of  $\tauolsint$. This may be because there is one stratum with a very small \revise{control} group size, size of ten, resulting in large finite-sample bias and variance for $\tauolsint$.

We also study the changes of the regression-adjusted estimators and $95\%$ confidence intervals when the number of covariates $K$ increases, from five to nine. Figure~\ref{fig:real_change_k} shows the results. The $\tauols$ produces stable point and interval estimates, with slightly decreasing confidence interval lengths. The $\tauolsint$ becomes worse and worse, which seems not to be consistent when the number of covariates is eight or nine.

\section{Discussion}
We study the theoretical advantage of regression adjustment in stratified randomized experiments. We assume that the number of covariates is fixed, not depending on the sample size. In practice, however, the number of covariates can converge to infinity at a rate of the sample size. In other words, the covariates are high-dimensional. For example, in a clinical trial, each patient's demographic and genetic information may be recorded. Performing high-dimensional regression adjustment and investigating its asymptotic properties in stratified randomized experiments are still challenging. The most difficult technical aspect  is to establish an appropriate concentration inequality for the mean of a stratified random sample, extending the Massart concentration inequality for the mean of a simple random sample used in completely randomized experiments \citep{bloniarz2015lasso,liu2018,Yue2018}. 

This paper focuses on using covariates adjustment to infer the average treatment effect for a binary treatment in stratified randomized experiments. It would be interesting to  extend the results to more complicated settings, such as covariates adjustment for multiple value treatment \citep{freedman2008regression_a},  for binary outcomes using logistic regression \citep{freedman2008randomization, zhang2008binary, moore2009}, and covariate adjustment in experiments when there is noncompliance \citep{imbens1994identification, angrist1995two, angrist1996identification}. Moreover, we proposed Neyman-type conservative variance estimators to construct large-sample conservative confidence intervals for the average treatment effect. It would be also interesting to explore other variance estimators under the randomized-based inference framework for stratified and sequentially randomized experiments, such as the Huber--White robust variance estimator for linear models.

\section*{Acknowledgement}
The authors thank Dr. Ke Zhu for his edits and suggestions that have helped clarify the text. Dr. Hanzhong Liu's research is supported by the National Natural Science Foundation of China (Grant No. 11701316). Dr. Yuehan Yang's research is partially supported by the National Natural Science Foundation of China (Grant No. 11671059).

\bibliographystyle{apalike}
\bibliography{causal}

\appendix

\section{Proof of main theorems}

\subsection*{Proof of Theorem~\ref{thm:fclt}}

\begin{proof}
Our proof relies on the following finite population central limit theorem \revise{proved} by \citet{bickel1984asymptotic}, who assumed a slightly weaker Lindeberg--Feller condition. Let $ v_i^2 $ be the ``variance weight": $v_i^2  = c_i^2 \S^2_{iy}(1)  \nci  / \{ \nti \ni  \sigma^2_{y}(1) \} = \var \{ c_i \hat y_{i\cdot}(1) / \sigma_{y}(1) \} $, and $\rho_i$ be the ``the effective sample size": $\rho_i = \nti \nci /\ni  $.
\begin{condition}(Lindeberg--Feller condition)
\label{cond:LF}
For any $\epsilon > 0$, suppose that
\begin{equation}
\label{eqn:LF}
\sumi  \frac{1} { \ni - 1 } \sum\limits_{j: |  v_i \{ y_{ij}(1) - y_{i\cdot}(1)  \} |  > \epsilon \S_{iy}(1)  \rho_i^{1/2}  } \frac{ v_i^2 \{ y_{ij}(1) - y_{i\cdot}(1) \}^2 }{ \S^2_{iy}(1) } \rightarrow 0. \nonumber
\end{equation}
\end{condition}
\begin{theorem}(Bickel and Freedman, 1984)
If Condition \ref{cond:LF} holds, then $ \{ \hat y_{\cdot \cdot}(1) - y_{\cdot \cdot } (1) \} / \sigma_{y}(1)$  converges in distribution to $\mathcal{N}(0,1)$. Furthermore, if for each stratum $i$, $2 \leq \nti \leq \ni - 1$, then $\hat \sigma^2_{y}(1) /\sigma^2_{y}(1)$ converges to one in probability.
\end{theorem}
What we have to show is that the condition $m_{1N} / N \rightarrow 0$ implies the Lindeberg--Feller condition. If
\begin{equation}
\label{eqn:maxij}
\lim\limits_{N \rightarrow \infty} \maxi \maxj \frac{v_i^2 \{ y_{ij}(1) - y_{i\cdot}(1) \}^2}{\S^2_{iy}(1) \rho_i} = 0,
\end{equation}
then for any $\epsilon > 0 $, there exists $N_{\epsilon}$ such that when $N \geq N_{\epsilon}$, $v_i^2 \{ y_{ij}(1) - y_{i\cdot}(1) \}^2 \leq  \epsilon^2 \S^2_{iy}(1) \rho_i$ for all $i= 1,\dots, B$ and $j = 1,\dots, \ni$. Thus, the second summation in the Lindeberg--Feller condition is zero, and hence, this condition holds. Recall that $c_i = \ni /N$ and $p_i = \nti / \ni$. By definite and simple calculation, we have
\begin{eqnarray}
\frac{v_i^2 \{ y_{ij}(1) - y_{i\cdot}(1) \}^2}{\S^2_{iy}(1) \rho_i}  & = & \frac{c_i^2 \S^2_{iy}(1)  \nci / (\nti\ni) \{ y_{ij}(1) - y_{i\cdot}(1) \}^2 }{ \S^2_{iy}(1) \big\{ \nti   \nci  / \ni \big\}  \sumi c_i^2 \S^2_{iy}(1)  \nci / (\nti \ni)     } \nonumber \\
& \leq &  \frac{1}{ p_i^2}   \frac{ \maxj \{  y_{ij}(1) - y_{i\cdot}(1) \}^2 }{N \sumi c_i \S^2_{iy}(1)  \nci /\nti  }  \leq  m_{1N}/N. \nonumber
\end{eqnarray}
Therefore, the condition $m_{1N}/N \rightarrow 0$ as $N \rightarrow \infty$ implies that the condition \eqref{eqn:maxij} holds, which further implies that the Lindeberg--Feller condition holds.

\end{proof}

\subsection*{Proof of Theorem~\ref{thm:clt}}

\begin{proof}
We write $\hat \tau_{i \cdot} - \tau_{i \cdot }$ as the average of a stratified random sample:
\begin{eqnarray}
\hat \tau_{i \cdot}  - \tau_{i \cdot } & = & \sumj \Big\{ \frac{Z_{ij}y_{ij}(1)}{\nti}  -   \frac{(1 - Z_{ij}) y_{ij}(0)}{\nci}   \Big\} - \{ y_{i\cdot}(1) - y_{i \cdot }(0) \} \nonumber \\
& = & \sumj \Big[ \frac{Z_{ij} \{ y_{ij}(1) - y_{i \cdot }(1) \} }{\nti}  -   \frac{(1 - Z_{ij}) \{ y_{ij}(0)  - y_{i \cdot }(0) \} }{\nci}   \Big]  \nonumber \\
& = &\frac{1}{ \nti } \sumj Z_{ij} \Big[ y_{ij}(1) - y_{i \cdot }(1)  + \frac{ \nti }{\nci}  \{ y_{ij}(0)  - y_{i \cdot }(0) \} \Big] , \nonumber
\end{eqnarray}
where the second equality is because $\sumj Z_{ij}/\nti = \sumj (1 - Z_{ij} )/ \nci = 1$ and the last equality is due to $\sumj \{y_{ij}(0)  - y_{i \cdot }(0) \} = 0 $. Let $a_{ij} = y_{ij}(1) - y_{i \cdot }(1)  + ( \nti /\nci )  \{ y_{ij}(0)  - y_{i \cdot }(0) \} $, then $ \hattaui - \tau_{i \cdot}  $ is the stratum-specific sample mean of a stratified random sample drawn from the population $\Pi_{a} = \{ a_{ij}: \   i = 1,\dots, B; \  j= 1,\dots, \ni \}$. Thus, $\tauunadj  - \tau = \sumi c_i ( \hattaui - \taui )$ is the weighted sample mean. We can apply Theorem~\ref{thm:fclt} to a stratified random sample drawn from $\Pi_{a}$ with the stratum-specific population mean $a_{i \cdot } = \sumj a_{ij}/ \ni =  0$. It is easy to show that the stratum-specific variance of $a_{ij}$ is
\begin{eqnarray}
\S_{ia}^2 & = & \frac{1 }{\ni - 1} \sumj  ( a_{ij} - a_{i\cdot} )^2 \nonumber \\
& = &  \frac{1 }{\ni - 1} \sumj  \Big[  y_{ij}(1) - y_{i \cdot }(1)  + \frac{ \nti }{\nci}  \{ y_{ij}(0)  - y_{i \cdot }(0) \}  \Big]^2 \nonumber \\
& = &  \S^2_{iy}(1) +  \frac{ \nti^2 }{ \nci^2 }  \S^2_{iy}(0) +  \frac{ \nti }{\nci}  \frac{2 }{\ni - 1} \sumj  \{  y_{ij}(1) - y_{i \cdot }(1) \}   \{  y_{ij}(0) - y_{i \cdot }(0) \}.  \nonumber
\end{eqnarray}
The above covariance term can be further expressed as
\begin{eqnarray}
& & - \frac{2 }{\ni - 1} \sumj  \big\{  y_{ij}(1) - y_{i \cdot }(1) \big\}   \big\{  y_{ij}(0) - y_{i \cdot }(0) \big\}  \nonumber \\
 & = & \frac{1 }{\ni - 1}  \sumj \Big[ \{ y_{ij}(1) - y_{i \cdot }(1) \} - \{ y_{ij}(0) - y_{i \cdot }(0) \} \Big]^2  -   \frac{1}{\ni - 1} \sumj \big\{  y_{ij}(1) - y_{i \cdot }(1) \big\}^2 \nonumber \\
 && -  \frac{1}{\ni - 1} \sumj \big\{  y_{ij}(0) - y_{i \cdot }(0) \big\}^2    \nonumber \\
 & = &  \frac{1 }{\ni - 1}  \sumj (  \tau_{ij} - \tau_{i \cdot } )^2 - \S^2_{iy}(1) - \S^2_{iy}(0) \nonumber \\
 & = & \S^2_{ i \tau}  - \S^2_{iy}(1) - \S^2_{iy}(0).  \nonumber
\end{eqnarray}
Since $\nti + \nci = \ni$, then
\begin{eqnarray}
\S_{ia}^2 & = &  \S^2_{iy}(1) +  \frac{ \nti^2 }{ \nci^2 }  \S^2_{iy}(0) + \frac{ \nti }{\nci}  \S^2_{iy}(1) + \frac{ \nti }{\nci}  \S^2_{iy}(0) - \frac{ \nti }{\nci}  \S^2_{ i \tau} = \frac{ \ni }{\nci} \S^2_{iy}(1) + \frac{ \nti \ni }{ \nci^2 }  \S^2_{iy}(0) -  \frac{ \nti }{\nci}  \S^2_{i \tau}. \nonumber
\end{eqnarray}
Since the randomization is independent across strata, then, the variance of $\tauunadj - \tau$ is
\begin{eqnarray}
\sigma_{\unadj}^2 = \sumi c_i^2 \var( \hat \tau_{i \cdot}  - \tau_{i \cdot }  ) =  \sumi c_i^2 \S_{ia}^2 \frac{\nci  }{\nti \ni } = \sumi c_i^2 \Big\{ \frac{\S^2_{iy}(1)}{\nti }  + \frac{\S^2_{iy}(0)}{\nci } - \frac{\S^2_{i \tau}}{\ni }   \Big\}. \nonumber
\end{eqnarray}
Moreover,
\begin{eqnarray}
\sumi c_i \S^2_{ia} \frac{\nci}{\nti} = \sumi c_i \Big\{  \frac{\S^2_{iy}(1)}{p_i} + \frac{\S^2_{iy}(0)}{1 - p_i} - \S^2_{i \tau}  \Big\} \nonumber
\end{eqnarray}
has a finite \revise{and positive} limit according to Condition \ref{cond:moments_outcomes}. Together with Condition \ref{cond:prop}, $p_i \rightarrow p_{i,\infty} \in (C, 1 - C)$ uniformly, the condition for asymptotic normality of a stratified random sample drawn from $\Pi_a$  is equivalent to
$$
\frac{1}{N} \maxi \maxj ( a_{ij} - a_{i\cdot} )^2  \rightarrow 0,
$$
which is implied by Condition \ref{cond:prop} and \ref{cond:max_dist_y}.

Next, we prove the asymptotical conservativeness of the variance estimator $ \sunadj^2$. The rescaled treated units $\{ (\ni / \nci)^{1/2} y_{ij}(1):  \  Z_{ij} = 1; \ i = 1,\dots, B; \ j = 1,\dots, \ni \}$ can be \revise{treated} as a stratified random sample drawn from the rescaled population $ \{  (\ni / \nci)^{1/2}  y_{ij}(1): \ i = 1,\dots, B; \ j = 1,\dots, \ni\}$, then \revise{applying} the second statement of Theorem \ref{thm:fclt}, we have
$
\{ \sumi c_i^2  \s^2_{iy}(1) / \nti  \} /  \{ \sumi c_i^2  \S^2_{iy}(1) / \nti  \}
$
tends to one in probability.  Similarly, \revise{applying} Theorem \ref{thm:fclt} to the rescaled control units $\{ (\ni / \nti)^{1/2} y_{ij}(0):  \  Z_{ij} = 0; \ i = 1,\dots, B; \ j = 1,\dots, \ni \}$, we have
$
\{ \sumi c_i^2  \s^2_{iy}(0) / \nci  \} /  \{ \sumi c_i^2  \S^2_{iy}(0) / \nci  \}
$
tends to one in probability. Therefore, the ratio
$$
\sumi c_i^2 \{ \s^2_{iy}(1) / \nti + \s^2_{iy}(0) / \nci  \} \Big/  \sumi c_i^2 \{ \S^2_{iy}(1)/ \nti +  \S^2_{iy}(0) / \nci  \}
$$
tends to one in probability. Under Condition \ref{cond:moments_outcomes}, the variance $\Sunadj^2$ is of order $1/N$ because
$$
\Sunadj^2 = \sumi c_i \frac{\ni}{N} \Big\{  \frac{\S^2_{iy}(1)}{ \nti } +   \frac{\S^2_{iy}(0)}{ \nci } - \frac{ \S_{\tau i}^2 }{ \ni }  \Big\} = \frac{1}{N} \sumi c_i \Big\{  \frac{\S^2_{iy}(1)}{ p_i } +   \frac{\S^2_{iy}(0)}{ 1 - p_i} -  \S_{i \tau}^2  \Big\}.
$$
Therefore, $N \sunadj^2$ converges in probability to $  \sumi c_i \{ \S^2_{iy}(1)  / p_i  + \S^2_{iy}(0)  / (1 - p_i )  \} $, and the difference $ N (  \sunadj^2 - \Sunadj^2  )  $ tends to \revise{$ \sumi c_i  \S_{i \tau}^2   \geq 0$.}

\end{proof}


\subsection*{Some useful lemmas}

\revise{To} obtain the asymptotic normality of $\tauols$, we require the following lemmas regarding the consistency of regression-adjusted coefficients $\hat \beta_1$, $\hat \beta_0$ and the asymptotic normality of stratified random samples from the covariates.

Let $\Pi_b = \{b_{ij}: \  i=1,\dots,B; \ j=1,\dots, \ni \}$ and $\Pi_d = \{d_{ij}: \  i=1,\dots,B; \ j=1,\dots, \ni \}$ be two series of fixed population quantities satisfying the following conditions:
\begin{condition}
\label{cond:lem1}
Let $\hat b_{i\cdot}$ and $\hat d_{i\cdot} $ be the  sample means of treated units in stratum $i$. As $N \rightarrow \infty$,

(a) the maximum squared distances satisfy
\[
\frac{1}{N} \maxi \maxj ( b_{ij} - b_{i\cdot} )^2  \rightarrow 0, \quad \frac{1}{N} \maxi \maxj ( d_{ij} - d_{i\cdot} )^2  \rightarrow 0;
\]

(b) the weighted variances,  $\sumi c_i  \S^2_{ib} $, $\sumi c_i \S^2_{id} $, and the weighted covariance,  $\sumi c_i  \S_{ibd}$, converge to finite limits (positive for variances), where $\S^2_{ib}$, $\S^2_{id}$ and $\S_{ibd}$ are the population variances of $b$, $d$ and the population covariance of $b$ and $d$ in stratum $i$, respectively.
\end{condition}

In our notations, $\S_{ibb} = \S^2_{ib}$ and  $\S_{idd} = \S^2_{id}$. Let  $\s_{ibd}$ be the sample covariance of $b$ and $d$ in stratum $i$:
\[
 \quad \s_{ibd} = \frac{1} { \nti -1 } \sumj Z_{ij} ( b_{ij} - \hat b_{i\cdot} ) ( d_{ij} -  \hat d_{i\cdot} ) .
\]
In what follows, $b_{ij}$ and $d_{ij}$ can be  $y_{ij}(1)$, $y_{ij}(0)$ and $\X_{ijk}$, for $k=1,\dots,K$.

\begin{lemma}
\label{lem:lem1}
Under Conditions  \ref{cond:prop} and \ref{cond:lem1}, \revise{and suppose that $\nti \geq 2$}, $\sumi c_i \s_{ibd} - \sumi c_i \S_{ibd} $ converges to zero in probability.
\end{lemma}

\begin{remark}
\revise{For the case with a few large strata, Lemma~\ref{lem:lem1} follows directly from Proposition 3 of \citet{fcltxlpd2016}. We will show that it also holds when there are a large number of small strata.}
\end{remark}

\begin{lemma}
\label{lem:consistent_beta}
Under Conditions \ref{cond:prop}--\ref{cond:moments_x}, we have $\hat \beta_1 - \beta_1$ and $\hat \beta_0 - \beta_0$ converge to zero in probability.
\end{lemma}

The proofs of Lemma~\ref{lem:lem1} and Lemma~\ref{lem:consistent_beta} are left to the last but second section of this supplementary material. Let
\begin{eqnarray}
\S^2_{ \X \X k} = \sumi c_i^2 \frac{1}{\ni - 1}  \sumj ( \x_{ijk} - \x_{i \cdot k}  )^2 \Big( \frac{1}{\nti} - \frac{1}{\ni} \Big)
 = \frac{1}{N} \sumi c_i  \frac{1}{\ni - 1}  \sumj ( \x_{ijk} - \x_{i \cdot k}  )^2 \frac{\nci}{\nti}, \nonumber
\end{eqnarray}
which is of order $1/N$ according to Conditions \ref{cond:prop} and \ref{cond:moments_x}. Applying the finite population central limit \revise{theorem} for stratified random sampling (Theorem~\ref{thm:fclt}) to each covariate $\X_{ijk}$ ($k= 1, \dots, K$), we can obtain the following lemma.
\begin{lemma}
\label{lem:normality_x}
Under Conditions \ref{cond:prop}--\ref{cond:moments_x}, for each covariate $\X_{ijk}$ ($k=1,\dots,K)$,  $\sumi c_i \{  \hat \X_{i \cdot k}(1) - \X_{i \cdot k}  \} / \S_{\X \X k}$  converges in distribution to $\mathcal{N}(0,1)$.
\end{lemma}

\subsection*{Proof of Theorem~\ref{thm:clt_ols}}

\begin{proof}

We first prove the asymptotic normality of $\tauols$. \revise{Recalling the projections of potential outcomes,}
\begin{equation}
\label{eqn:decomp_appendix}
y_{ij}(1) = y_{i\cdot}(1) + \big( \X_{ij} - \X_{i\cdot} \big)^\T \beta_{1} + \v_{ij}(1), \quad y_{ij}(0) = y_{i\cdot}(0) + \big( \X_{ij} - \X_{i\cdot} \big)^\T \beta_{0} + \v_{ij}(0).
\end{equation}
It is easy to verify that $\v_{i\cdot}(1) = \v_{ i\cdot } (0)= 0$. If we define $h_1 = \hat \beta_1 - \beta_1$, $h_0 = \hat \beta_0 - \beta_0$, by substitution, we have
\begin{eqnarray}
\tauiols & = &   \Big[  \hat y_{i\cdot}(1) - \big\{ \hat X_{i\cdot}(1)  - X_{i\cdot} \big\}^\T \hat \beta_{1} \Big]  - \Big[ \hat y_{i\cdot}(0) - \big\{ \hat X_{i\cdot}(0) - X_{i\cdot} \big\}^\T \hat \beta_{0} \Big] \nonumber \\
& = & \Big[ \big \{ \hat \v_{i\cdot}(1) - \hat \v_{i\cdot}(0) \big \} + \tau_{i \cdot}  \Big]  - \Big[  \big\{ \hat X_{i\cdot}(1)  -  X_{i\cdot} \big\}^\T h_1 - \big\{ \hat X_{i\cdot}(0) - X_{i\cdot} \big\}^\T h_0  \Big]. \nonumber
\end{eqnarray}
Thus,
\begin{eqnarray}
\label{eqn:two-terms}
\frac{ \tauols - \tau }{ \sigma_{\ols} } =  \frac{ \sumi c_i \big\{ \hat \v_{i\cdot}(1) - \hat \v_{i\cdot}(0) \big\} }{ \sigma_{\ols} }    - \frac{ \sumi c_i  \Big[  \big\{ \hat X_{i\cdot}(1)  -  X_{i\cdot} \big\}^\T h_1 - \big\{ \hat X_{i\cdot}(0) - X_{i\cdot} \big\}^\T h_0  \Big] } { \sigma_{\ols} } .
\end{eqnarray}
We analyze the two terms on the right hand side separately, showing that the first term is asymptotically normal and the second term converges to zero in probability.

The first term  is the stratified difference-in-means estimator divided by its standard deviation  for the potential outcomes $\v_{ij}(1)$ and $\v_{ij}(0)$, whose population average treatment effect is $\sumi c_i \{ \v_{ i\cdot } (1)- \v_{ i\cdot } (0)\} = 0$. By the definition of projection coefficient \revise{vector} $\beta_1 = (\beta_{11}, \dots, \beta_{1K})^\T$, we have
\[
\beta_1 = \argmin\limits_{\beta} \sumi \frac{c_i}{\ni - 1} \sumj \Big[ y_{ij}(1) - y_{i\cdot}(1) - \big\{ \X_{ij} - \X_{i\cdot}  \big\}^\T \beta \Big]^2 =   \S_{\X \X}^{-1}\S_{\X y}(1),
\]
which has a finite limit (a $K$-dimensional vector) as $N \rightarrow \infty$ according to Condition~\ref{cond:moments_x}. Thus, the maximum squared distance for $\v_{ij}(1)$ satisfies
\begin{eqnarray}
& & \maxi \maxj \big\{ \v_{ij}(1) - \v_{i\cdot}(1) \big \}^2  \nonumber \\
& \leq & 2 \maxi \maxj \big\{ y_{ij}(1) - y_{i\cdot}(1) \big\}^2 + 2 K^2 \max\limits_{k=1,\dots,K} \maxi \maxj ( \x_{ijk} - \x_{i\cdot k} )^2  \max\limits_{k=1,\dots,K} | \beta_{1k}|^2 \nonumber.
\end{eqnarray}

Therefore, $\max_{i=1,\dots,B} \max_{j=1,\dots,\ni } \{ \v_{ij}(1) - \v_{i\cdot}(1) \}^2 /N$ tends to zero due to Conditions~\ref{cond:max_dist_y} and \ref{cond:max_dist_x} (the maximum squared distance for $y_{ij}(1)$ and $\X_{ijk}$). Similarly,  $\max_{i=1,\dots,B} \max_{j=1,\dots,\ni } \{ \v_{ij}(0) - \v_{i\cdot}(0) \}^2 /N$ tends to zero. Therefore, the maximum squared distance \revise{conditions} for potential outcomes $\v_{ij}(1)$ and $\v_{ij}(0)$ \revise{hold}. Applying the finite population central limit theorem for stratified random sampling (Thereom~\ref{thm:clt}) to $\v_{ij}(1)$ and $\v_{ij}(0)$, we obtain the asymptotic normality of the first term, i.e., $\sumi c_i \{ \hat \v_{i\cdot}(1) - \hat \v_{i\cdot}(0)\} / \sigma_{\ols} $ converges  in distribution to $\mathcal{N}(0,1)$.

For the second term in \eqref{eqn:two-terms}, Lemma~\ref{lem:consistent_beta} implies that both $h_1$ and $h_0$ converge to zero in probability and Lemma~\ref{lem:normality_x} implies that each coordinate, $ \sumi c_i \{ \hat \X_{i\cdot k}(1)  -  \X_{i\cdot k}  \} /  \S_{\X \X k}$, converges  in distribution to $\mathcal{N}(0,1)$. By Condition~\ref{cond:variance_e}, the variance $\sigma^2_{\ols}$ is of order $1/N$ (Similar to the argument for $\Sunadj^2$). Since both $    \S^2_{\X \X k}$ and $ \sigma^2_{\ols}$ are of order $1/N$, the second term converges to zero in probability.

Second, we compare the asymptotic variance of $\surd{N} \tauols$ and $\surd{N} \tauunadj$. For stratum $i$,
\begin{eqnarray}
\label{eqn:S1iy}
\S^2_{iy}(1) & = & \frac{1}{\ni - 1 } \sumj \big\{ y_{ij}(1) - y_{i\cdot}(1) \big\}^2 \nonumber \\
& = &  \frac{1}{\ni - 1 }  \sumj \Big\{ \big(  \X_{ij} - \X_{i\cdot}  \big)^\T \beta_1 \Big\}^2 + \frac{1}{\ni - 1} \sumj \v^2_{ij}(1) + \frac{2}{\ni - 1} \sumj  \v_{ij}(1) \big(  \X_{ij} - \X_{i\cdot}  \big)^\T \beta_1 \nonumber \\
& = & \beta_1^\T \S_{i\X\X} \beta_1 + \S^2_{i\v}(1) + 2 \S_{i \X \v}^\T(1) \beta_1.
\end{eqnarray}
Similarly,
\begin{eqnarray}
\label{eqn:S0iy}
\S^2_{iy}(0)  =   \beta_0^\T \S_{i\X\X} \beta_0 + \S^2_{i\v}(0) + 2 \S_{i \X \v}^\T(0) \beta_0,
\end{eqnarray}
\begin{eqnarray}
\label{eqn:Sitau}
\S^2_{i \tau}  =   (\beta_1 - \beta_0 )^\T \S_{i\X\X} (\beta_1 - \beta_0) + \S^2_{i\v} + 2 \big\{ \S_{i \X \v}(1) - \S_{i\X \v}(0) \big\}^\T \big( \beta_1 - \beta_0 \big) .
\end{eqnarray}
Thus, simple calculus gives
\begin{eqnarray}
&& \Big\{ \frac{ \S^2_{iy}(1)  }{\nti} + \frac{ \S^2_{iy}(0)  }{\nci} - \frac{\S^2_{i \tau}  }{\ni}  \Big\} - \Big\{ \frac{ \S^2_{i \v}(1)  }{\nti} + \frac{ \S^2_{i \v}(0)  }{\nci} - \frac{\S^2_{i \v}  }{\ni}  \Big\} \nonumber \\
& = & \frac{1}{\ni} \frac{1 }{ p_i ( 1 - p_i ) }   \big( \beta^i \big)^\T \S_{i \X \X} \big( \beta^i \big) + \frac{2}{\ni} \frac{1}{p_i} \S_{i\X \v}^\T(1) \beta^i +  \frac{2}{\ni} \frac{1}{1 - p_i} \S_{i\X \v}^\T(0) \beta^i \nonumber \\
& = & \frac{1}{\ni} \Delta_i^2 + \frac{2}{\ni} \frac{1}{p_i} \S_{i\X \v}^\T(1) \beta^i +  \frac{2}{\ni} \frac{1}{1 - p_i} \S_{i\X \v}^\T(0) \beta^i, \nonumber
\end{eqnarray}
where $p_i = \nti/\ni  \rightarrow p$, $\beta^i = ( 1 - p_i) \beta_1 + p_i \beta_0 \rightarrow \lim_{N \rightarrow \infty} (1 - p) \beta_1 + p \beta_ 0$, and
\[ \Delta_i = 1/\{p_i (1 - p_i ) \}    \big( \beta^i \big)^\T \S_{i \X \X} \big( \beta^i \big) . \]
Therefore, the difference between the asymptotic variance of $\surd{N} \tauunadj$ and $\surd{N} \tauols$ is
\begin{eqnarray}
N \sigma^2_{\unadj} - N \sigma^2_{\ols} & = & N \sumi c_i^2 \frac{1}{\ni} \Delta_i^2  + N \sumi c_i^2  \frac{2}{\ni} \frac{1}{p_i} \S_{i\X \v}^\T(1) \beta^i + N \sumi c_i^2  \frac{2}{\ni} \frac{1}{1 - p_i} \S_{i\X \v}^\T(0) \beta^i \nonumber \\
& = & \sumi c_i \Delta_i^2  +  2 \sumi c_i  \frac{1}{p_i} \S_{i\X \v}^\T(1) \beta^i + 2 \sumi c_i \frac{1}{1 - p_i} \S_{i\X \v}^\T(0) \beta^i . \nonumber
\end{eqnarray}
It is enough to show that as $N \rightarrow \infty$,
\begin{equation}
\label{eqn:small_order_term}
\sumi c_i  \frac{1}{p_i} \S_{i\X \v}^\T(1) \beta^i \rightarrow 0 , \quad \sumi c_i \frac{1}{1 - p_i} \S_{i\X \v}^\T(0) \beta^i \rightarrow 0.
\end{equation}
By the definition of projection (or the definition of projection coefficient \revise{vector} $\beta_1$), the projection errors $\v_{ij}(1)$ are orthogonal to the covariates $\X_{ij}$ in the following sense:
\[
\sumi c_i \frac{1}{\ni -1} \sumj ( \X_{ij} - \X_{i \cdot} )  \v_{ij}(1)  = \sumi c_i \S_{i\X \v}(1) = 0.
\]
On the other hand,
\[
 \S_{i\X \v}(1) = \S_{i \X y}(1) - \S_{i\X \X} \beta_1.
\]
Therefore, if let $\beta_{\infty} =  \lim_{N \rightarrow \infty} \beta^i = \lim_{N \rightarrow \infty} ( 1 - p ) \beta_1 + p \beta_0$, then
\begin{eqnarray}
\sumi c_i  \frac{1}{p_i} \S_{i\X \v}^\T(1) \beta^i  & = &  \sumi c_i  \frac{1}{p_i} \S_{i\X \v}^\T(1) \beta^i - \sumi c_i \S_{i\X \v}^\T(1)  \frac{\beta_{\infty}}{p}  \nonumber \\
 & = & \sumi c_i \S_{i\X \v}^\T(1) \Big(   \frac{1}{p_i} \beta^i - \frac{\beta_{\infty}}{p}  \Big) \nonumber \\
 & \leq & \sumi c_i | \S_{i \X y}(1)  |^\T  \Big |  \Big(   \frac{1}{p_i} \beta^i - \frac{\beta_{\infty}}{p}  \Big) \Big | +  | \beta_1 |^\T \sumi c_i | \S_{i \X \X}  |^\T  \Big |  \Big(   \frac{1}{p_i} \beta^i - \frac{\beta_{\infty}}{p}  \Big) \Big |, \nonumber
\end{eqnarray}
which tends to zero because $\beta^i / p_i - \beta_{\infty} / p $ tends to zero \revise{uniformly for $i=1,\dots,B$}, and the weighted absolute covariances $\sumi c_i | \S_{i \X y}(1)  |$, $\sumi c_i | \S_{i \X \X}  | $ tend to finite limits (Condition~\ref{cond:moments_x}). Similarly,  $ \sumi c_i  \S_{i\X \v}^\T(0) \beta^i /(1 - p_i ) $ tends to zero.

Finally, we prove the conservativeness of the variance estimator $\hat \sigma^2_{\ols}$. Let $\s_{i \X \X}(1)$, $  \s_{i\X \v}(1) $, and $\s_{i \X \X}(0)$, $\s_{i\X \v}(0)$ be the sample covariance matrices (or covariances) computed from the \revise{treatment and control groups in stratum $i$}, respectively. For example,
\begin{align*}
  &\s_{i \X \X}(1) = \frac{1}{\nti - 1}  \sumj Z_{ij}   \big\{ \X_{ij} -  \hat \X_{i\cdot}(1) \big\}   \big\{ \X_{ij} -  \hat \X_{i\cdot}(1) \big\}^\T,\\
  &\s_{i\X \v}(1) = \frac{1}{\nti - 1}  \sumj Z_{ij}  \big\{ \X_{ij} -  \hat \X_{i\cdot}(1) \big\} \big\{ \v_{ij}(1) - \hat \v_{i \cdot}(1) \big\}.
\end{align*}
By definition and \revise{projections} of potential outcomes  \eqref{eqn:decomp_appendix}, we have
\begin{eqnarray}
\hat \s^2_{i \v}(1) & = & \frac{1}{\nti - 1} \sumj Z_{ij}  \Big[ y_{ij}(1) - \hat y_{i\cdot}(1) - \big\{ \X_{ij} -  \hat \X_{i\cdot}(1) \big\}^\T \hat \beta_{1} \Big]^2  \nonumber \\
& = &  \frac{1}{\nti - 1} \sumj Z_{ij} \Big[ \big\{ \X_{ij} -  \hat \X_{i\cdot}(1) \big\}^\T \big(  \beta_1 - \hat \beta_1  \big)  + \v_{ij}(1) - \hat \v_{i \cdot}(1) \Big]^2 \nonumber \\
& = & \big(  \beta_1 - \hat \beta_1  \big)^\T \s_{i \X \X}(1) \big(  \beta_1 - \hat \beta_1  \big)  + \s^2_{i\v}(1) + \big(  \beta_1 - \hat \beta_1  \big)^\T \s_{i\X \v}(1). \nonumber
\end{eqnarray}
Thus,
\begin{eqnarray}
& & N  \sumi c_i^2 \frac{\hat \s^2_{i \v}(1) }{\nti} - N \sumi c_i^2 \frac{ \S^2_{i\v}(1) }{\nti}  \nonumber \\
& = &  \big(  \beta_1 - \hat \beta_1  \big)^\T \sumi c_i \frac{\ni}{\nti}  \s_{i \X \X}(1) \big(  \beta_1 - \hat \beta_1  \big) + \Big\{  \sumi c_i \frac{ \s_{i\v}^2(1)}{p_i }  - \sumi c_i \frac{ \S^2_{i\v}(1) }{p_i}  \Big\} \nonumber \\
& & +  \big(  \beta_1 - \hat \beta_1  \big)^\T \sumi c_i \s_{i\X \v}(1) \frac{\ni}{\nti}. \nonumber
\end{eqnarray}
Therefore,
\begin{eqnarray}
& & N \hat \sigma^2_{\ols} - \sumi  c_i \Big\{ \frac{\S^2_{i\v}(1) }{p_i}  +  \frac{\S^2_{i\v}(0) }{1 - p_i} \Big\} \nonumber \\
& = &  \big(  \beta_1 - \hat \beta_1  \big)^\T \sumi c_i \frac{\ni}{\nti}  \s_{i \X \X}(1) \big(  \beta_1 - \hat \beta_1  \big) +  \big(  \beta_0 - \hat \beta_0  \big)^\T \sumi c_i \frac{\ni}{\nci}  \s_{i \X \X}(0) \big(  \beta_0 - \hat \beta_0  \big)  \label{eqn:part1} \\
&& +  \Big[  \sumi c_i  \Big\{ \frac{ \s^2_{i\v}(1) }{p_i} +  \frac{\s^2_{i\v}(0) }{1 - p_i} \Big\}  - \sumi c_i  \Big\{ \frac{ \S^2_{i\v}(1) }{p_i } + \frac{ \S^2_{i\v}(0) }{1 - p_i}  \Big\}  \Big]  \label{eqn:part2}  \\
& & +  \big(  \beta_1 - \hat \beta_1  \big)^\T \sumi c_i \s_{i\X \v}(1) \frac{\ni}{\nti} +  \big(  \beta_0 - \hat \beta_0  \big)^\T \sumi c_i \s_{i\X \v}(0) \frac{\ni}{\nti}.  \label{eqn:part3}
\end{eqnarray}

Applying Lemma~\ref{lem:lem1} to $b_{ij} = \X_{ij k}$ and $d_{ij} = \X_{ijk'}$, we can obtain the element-wise convergence of $\sumi c_i  \s_{i \X \X}(1) $ to $\S_{\X\X} = \sumi c_i   \S_{i \X \X}$ in probability. Since $\ni / \nti = 1/p_i \rightarrow 1/p$ uniformly, we have
\[
 \sumi c_i \frac{\ni}{\nti}  \s_{i \X \X}(1)  - \frac{1}{p} \sumi c_i   \S_{i \X \X}
\]
tends to zero in probability.  Lemma~\ref{lem:consistent_beta} implies \revise{that} $\beta_1 - \hat \beta_1$ converges to zero in probability element-wise. Therefore, the first term in the summation, \eqref{eqn:part1}, tends to zero in probability. Theorem~\ref{thm:clt} implies that the second term in the summation, \eqref{eqn:part2}, tends to zero in probability. Similarly, for the third term, \eqref{eqn:part3}, applying Lemma~\ref{lem:lem1} to $b_{ij} = \X_{ij k}$ and $d_{ij} = \v_{ij}(1)$ or $\v_{ij}(0)$, together with Lemma~\ref{lem:consistent_beta}, it tends to zero in probability. Therefore, $N \hat \sigma^2_{\ols} $ converges in probability to the limit of $\sumi  c_i  \{  \S^2_{i\v}(1) / p_i  +  \S^2_{i\v}(0) / (1 - p_i) \}$.

Since $N \sigma^2_{\ols} = \sumi  c_i  \{  \S^2_{i\v}(1) / p_i  +  \S^2_{i\v}(0) / ( 1 - p_i) - \S^2_{i\v} \}$, the difference $N \hat \sigma^2_{\ols} - N \sigma^2_{\ols}$ converges in probability to the limit of $ \sumi c_i  \S^2_{i\v}$, which is greater than or equal to zero.

By Theorem~\ref{thm:clt}, $N \sunadj^2$ converges to the limit of $ \sumi  c_i  \{  \S^2_{i y}(1) / p_i  +  \S^2_{i y}(0) / (1 - p_i)  \}$, thus,  the difference $N \sunadj^2 - N \hat \sigma^2_{\ols}$ converges in probability to  the limit of
\[
\sumi c_i \Big\{  \frac{\S^2_{iy}(1)} {p_i } + \frac{\S^2_{iy}(0)}{ 1 - p_i }  \Big\}  - \sumi c_i \Big\{ \frac{\S^2_{i\v}(1) }{p_i} + \frac{\S^2_{i\v}(0) }{ 1- p_i }  \Big\}.
\]
We have shown in \eqref{eqn:S1iy} and \eqref{eqn:S0iy} that
\[
\S^2_{iy}(1)  =   \beta_1^\T \S_{i\X\X} \beta_1 + \S^2_{i\v}(1) + 2 \S_{i \X \v}^\T(1) \beta_1,
\]
\[
\S^2_{iy}(0)  =   \beta_0^\T \S_{i\X\X} \beta_0 + \S^2_{i\v}(0) + 2 \S_{i \X \v}^\T(0) \beta_0.
\]
Thus,
\begin{eqnarray}
&& \sumi c_i \Big\{  \frac{\S^2_{iy}(1)} {p_i } + \frac{\S^2_{iy}(0)}{ 1 - p_i }  \Big\}  - \sumi c_i \Big\{ \frac{\S^2_{i\v}(1) }{p_i} + \frac{\S^2_{i\v}(0) }{ 1- p_i }  \Big\} \nonumber \\
& = &\beta_1^\T  \sumi c_i \frac{1}{p_i}  \S_{i\X \X} \beta_1 +  \beta_0^\T \sumi c_i \frac{1}{1 - p_i }  \S_{i\X \X} \beta_0  + 2 \sumi c_i   \frac{1}{p_i} \S^\T_{i \X \v}(1) \beta_1 + 2 \sumi c_i   \frac{1}{ 1 - p_i} \S^\T_{i \X \v} (0)\beta_0.  \nonumber
\end{eqnarray}
Similar to \eqref{eqn:small_order_term}, we can show that
\[
\sumi c_i   \frac{1}{p_i} \S^\T_{i \X \v}(1) \beta_1 \rightarrow 0 , \quad \sumi c_i   \frac{1}{ 1 - p_i} \S^\T_{i \X \v}(0) \beta_0 \rightarrow 0,
\]
and
\[
\beta_1^\T  \sumi c_i \frac{1}{p_i}  \S_{i\X \X} \beta_1  \rightarrow   \lim_{N \rightarrow \infty} \frac{1}{p} \beta_1^\T  \S_{\X\X} \beta_1, \quad  \beta_0^\T \sumi c_i \frac{1}{1 - p_i }  \S_{i\X \X} \beta_0 \rightarrow   \lim_{N \rightarrow \infty}  \frac{1}{1-p} \beta_0^\T \S_{\X\X} \beta_0.
\]
Therefore, the difference $N \sunadj^2 - N \hat \sigma^2_{\ols}$ converges in probability to  the limit of
\[
 \frac{1}{p} \beta_1^\T  \S_{\X\X} \beta_1 + \frac{1}{1-p} \beta_0^\T \S_{\X\X} \beta_0.
\]

\end{proof}

\subsection*{Proof of Theorem~\ref{thm:clt_ols_int}}

\begin{proof}
We use the following asymptotic normality result for regression-adjusted average treatment effect estimator in completely randomized experiments. See Example 9 in \citet{fcltxlpd2016} and also Theorem 1 and Corollary 1.1 in \citet{lin2013} who assumed \revise{the stronger fourth moment conditions}. Let 
\[
\sigma^2_{i,\olsint}  = \S^2_{i\e}(1)/n_{1i}  +  \S^2_{i\e}(0) / n_{0i} - \S^2_{i\e} / \ni, \quad \hat \sigma^2_{i,\olsint}  = \hat \s^2_{i\e}(1)/n_{1i}  +  \hat \s^2_{i\e}(0) / n_{0i}, 
\]
\[ \hat \sigma^2_{i,\unadj} = \s^2_{iy}(1)/n_{1i}  +   \s^2_{iy}(0) / n_{0i}.
\]

\begin{proposition}(\citet{fcltxlpd2016})
If Conditions \ref{cond:prop}, \ref{cond:max_dist_y}, \ref{cond:max_dist_x} and \ref{cond:large_ni} hold, then for stratum $i$, $( \tauiolsint - \taui ) /  \sigma_{i,\olsint}   $  converges in distribution to $\mathcal{N}(0,1)$, as $\ni \rightarrow \infty$. The difference between the asymptotic \revise{variances} of $\surd{\ni} ( \tauiolsint - \taui ) $ and $\surd{\ni} ( \hat \tau_{i\cdot}  - \taui )$ tends to the limit of  $- \tilde \Delta_i^2$, where
\[
 \tilde \Delta_i^2  =  \frac{ 1 }{ p_i  ( 1 - p_i ) }   \big( \tilde \beta^i \big)^\T \S_{i\X \X} \big( \tilde \beta^i \big),  \quad \tilde \beta^i  = (1 - p_i ) \beta_{1i} + p_i \beta_{0i}.
 \]
Furthermore, the  estimator $\ni \hat \sigma^2_{i,\olsint} $ converges in probability to the limit of $\S^2_{i\e}(1)/p_i + \S^2_{i\e}(0) / (1 - p_i)$, which is no less than that of $\ni \sigma^2_{i,\olsint}$, and the difference is the limit of $\S^2_{i \e}  $. The difference $\ni  ( \hat \sigma^2_{i, \olsint}  -  \hat \sigma^2_{i, \unadj} )$ converges in probability to the limit of
\[
- \frac{1}{p_i}  \beta_{1i}^\T \S_{i\X\X} \beta_{1i}   - \frac{1}{1-p_i} \beta_{0i}^\T \S_{i\X\X} \beta_{0i}.
\]
\end{proposition}

Since $( \tauiols - \taui )  $ are independent across strata, and the number of strata $B$ is finite, Theorem \ref{thm:clt_ols_int} holds immediately.

\end{proof}

\subsection*{Proof of Corollary~\ref{col:var_ols_int}}

\begin{proof}

We have shown in the proof of Theorem~\ref{thm:clt_ols} that $N \sigma^2_{\ols} \leq N \Sunadj^2$. Now we show that $N \sigma^2_{\olsint}  \leq N \sigma^2_{\ols} $. According to the projections of potential outcomes,  we have
\[
y_{ij}(1) = y_{i\cdot}(1) + \big( \X_{ij} - \X_{i\cdot} \big)^\T \beta_{1} + \v_{ij}(1), \quad y_{ij}(1) = y_{i\cdot}(1) + \big( X_{ij} - X_{i\cdot} \big)^\T \beta_{1i} + \e_{ij}(1).
\]
Taking difference gives
\[
 \v_{ij}(1) =  \big( \X_{ij} - \X_{i\cdot} \big)^\T \big( \beta_{1i} -  \beta_{1} \big) + \e_{ij}(1).
\]
Thus,
\begin{eqnarray}
\label{eqn:S1ive}
\S^2_{i \v}(1) & = & \frac{1}{\ni - 1 } \sumj \big\{ \v_{ij}(1) - \v_{i\cdot}(1) \big\}^2 \nonumber \\
& = &  \frac{1}{\ni - 1 }  \sumj \Big\{ \big(  \X_{ij} - \X_{i\cdot}  \big)^\T  \big( \beta_{1i} -  \beta_{1} \big) \Big\}^2 + \frac{1}{\ni - 1} \sumj \e^2_{ij}(1) \nonumber \\
&& + \frac{2}{\ni - 1} \sumj  \e_{ij}(1) \big(  \X_{ij} - \X_{i\cdot}   \big)^\T  \big( \beta_{1i} -  \beta_{1} \big)  \nonumber \\
& = &  \frac{1}{\ni - 1 }  \sumj \Big\{ \big(  \X_{ij} - \X_{i\cdot}  \big)^\T  \big( \beta_{1i} -  \beta_{1} \big) \Big\}^2 + \S^2_{i\e}(1) \nonumber \\
& = &  \big(  \beta_{1i} -  \beta_{1} \big)^\T \S_{i\X\X}  \big( \beta_{1i} -  \beta_{1} \big)  +  \S^2_{i\e}(1),
\end{eqnarray}
where the last but second equality is due to the property of projection that the covariates $\X_{ij} - \X_{i\cdot} $ and the projection \revise{errors} $\e_{ij}(1)$ are orthogonal in stratum $i$, i.e., $ \sumj  \e_{ij}(1) (  \X_{ij} - \X_{i\cdot}   ) = 0$. Similarly,
\begin{eqnarray}
\label{eqn:S0ive}
\S^2_{i \v}(0) & =   \big( \beta_{0i} -  \beta_{0} \big)^\T \S_{i\X\X}  \big( \beta_{0i} -  \beta_{0} \big)  + \S^2_{i\e}(0),
\end{eqnarray}
\begin{eqnarray}
\label{eqn:Sive}
\S^2_{i \v} & =   \big( \beta_{1i} - \beta_1 - \beta_{0i} +  \beta_{0} \big)^\T   \S_{i\X\X}     \big( \beta_{1i} - \beta_1 - \beta_{0i} +  \beta_{0} \big) + \S^2_{i\e}.
\end{eqnarray}
Combining \eqref{eqn:S1ive} -- \eqref{eqn:Sive}, we have
\begin{eqnarray}
 N \sigma^2_{\ols}  - N \sigma^2_{\olsint}  & = & \sumi c_i \frac{ 1 }{ p_i }   \Big(  \beta_{1i} -  \beta_{1} \Big)^\T \S_{i\X\X}  \Big( \beta_{1i} -  \beta_{1} \Big)   \nonumber \\
 &  & +   \sumi c_i \frac{ 1 }{ 1 - p_i }  \Big( \beta_{0i} -  \beta_{0} \Big)^\T \S_{i\X\X}  \Big( \beta_{0i} -  \beta_{0} \Big) \nonumber \\
 && - \sumi c_i   \Big( \beta_{1i} - \beta_1 - \beta_{0i} +  \beta_{0} \Big)^\T   \S_{i\X\X}     \Big( \beta_{1i} - \beta_1 - \beta_{0i} +  \beta_{0} \Big)  \nonumber \\
 & = & \sumi c_i \frac{1}{ p_i ( 1 - p_i ) } \gamma_i^\T \S_{i\X\X} \gamma_i  \geq 0, \nonumber
\end{eqnarray}
where $\gamma_i = ( 1  - p_i ) ( \beta_{1i} - \beta_1 ) + p_i ( \beta_{0i} - \beta_0 )$.

Next, we show $N \hat \sigma^2_{\olsint} \leq N \hat \sigma^2_{\ols} \leq N \sunadj^2$ holds in probability. We have shown in Theorem~\ref{thm:clt_ols} that $N \hat \sigma^2_{\ols} - N \sunadj^2$ converges in probability to a limit no larger than zero, and $N \hat \sigma^2_{\ols}$ converges in probability to the limit of $  \sumi c_i \{ \S_{i \v}^2(1)  / p + \S_{i \v}^2(0)  / ( 1 - p )  \}  $ (in this corollary we assume \revise{that} the proportion of treated units is asymptotically the same across strata, that is, $p_{i,\infty} = p$). Theorem~\ref{thm:clt_ols_int} shows that $N \hat \sigma^2_{\olsint}$ converges in probability to the limit of $ \sumi c_i \{ \S_{i \e}^2(1)  / p + \S_{i \e}^2(0)  / ( 1 - p )  \} $. Therefore, it is enough to show that the limit of
\[
\sumi c_i \Big\{ \S_{i \v}^2(1)  / p + \S_{i \v}^2(0)  / ( 1 - p )  \Big\} -   \sumi c_i  \Big \{ \S_{i \e}^2(1)  / p + \S_{i \e}^2(0)  / ( 1 - p ) \Big \}
\]
is no less than zero. Combing \eqref{eqn:S1ive} and \eqref{eqn:S0ive}, we have
\begin{eqnarray}
& & \sumi c_i \Big \{ \S_{i \v}^2(1)  / p + \S_{i \v}^2(0)  / ( 1 - p ) \Big \} -   \sumi c_i \Big \{ \S_{i \e}^2(1)  / p + \S_{i \e}^2(0)  / ( 1 - p )  \Big \} \nonumber \\
& = & \frac{1}{p}  \sumi c_i   \Big(  \beta_{1i} -  \beta_{1} \Big)^\T \S_{i\X\X}  \Big( \beta_{1i} -  \beta_{1} \Big)   + \frac{1}{1 - p}   \sumi c_i   \Big( \beta_{0i} -  \beta_{0} \Big)^\T \S_{i\X\X}  \Big( \beta_{0i} -  \beta_{0} \Big) \geq 0, \nonumber
\end{eqnarray}
since the covariance matrix $\S_{i\X\X} $ is positive-definite. 

\end{proof}

\section*{Proof of lemmas}
\label{app:proof_lemma}

\subsection*{Proof of lemma~\ref{lem:lem1}}

\begin{proof}
 Applying standard sampling theory (see for example \citet{cochran1977}), it is easy to show that the mean of $\s_{ibd}$ is $\S_{ibd}$, i.e., $E ( \s_{ibd} ) = \S_{ibd}$. Therefore,
 \[
 E \Big( \sumi c_i \s_{ibd} - \sumi c_i \S_{ibd} \Big) = 0 .
 \]
Using Markov inequality, it is enough for Lemma~\ref{lem:lem1} to show that the variance $\var (\sumi c_i \s_{ibd}) $ tends to zero. Since completely randomized experiment is  conducted independently across strata, we have $ \var (\sumi c_i \s_{ibd})  = \sumi c_i^2 \var( \s_{ibd}) $. Next, we bound the stratum-specific variance $\var( \s_{ibd})$. Since
\begin{eqnarray}
\s_{ibd} & =  & \frac{1} { \nti -1 } \sumj Z_{ij} ( b_{ij} - \hat b_{i\cdot} ) ( d_{ij} -  \hat d_{i\cdot} ) \nonumber \\
&  =  & \frac{1} { \nti -1 } \sumj Z_{ij} ( b_{ij} - b_{i\cdot} ) ( d_{ij} -  d_{i\cdot} )  - \frac{ \nti }{ \nti - 1} ( \hat b_{i\cdot } - b_{i\cdot} ) ( \hat d_{i \cdot} -  d_{i\cdot} ) , \nonumber
\end{eqnarray}
then
\begin{eqnarray}
\var( \s_{ibd}  ) & \leq  &  2 \var \Big(   \frac{1} { \nti -1 } \sumj Z_{ij} ( b_{ij} - b_{i\cdot} ) ( d_{ij} -  d_{i\cdot} )   \Big) + 2 \var \Big(  \frac{ \nti }{ \nti - 1} ( \hat b_{i\cdot } - b_{i\cdot} ) ( \hat d_{i \cdot} -  d_{i\cdot} )  \Big). \nonumber
\end{eqnarray}

The first term is bounded as follows:
\begin{eqnarray}
\label{eqn:first_term}
&&  \var \Big(   \frac{1} { \nti -1 } \sumj Z_{ij} ( b_{ij} - b_{i\cdot} ) ( d_{ij} -  d_{i\cdot} )   \Big)  \nonumber \\
& = &  \frac{ \nti^2 } { ( \nti -1 )^2 } \var \Big (  \frac{1}{\nti}  \sumj Z_{ij} ( b_{ij} - b_{i\cdot} ) ( d_{ij} -  d_{i\cdot} )   \Big) \nonumber \\
& \leq &  \frac{ \nti^2 } { ( \nti -1 )^2 } \frac{1}{ \ni - 1 }  \sumj \Big\{ ( b_{ij} - b_{i\cdot} ) ( d_{ij} -  d_{i\cdot} )   \Big\}^2 \Big( \frac{1}{\nti} - \frac{1}{\ni} \Big)  \nonumber \\
& \leq & \maxj  ( b_{ij} - b_{i\cdot} )^2 \frac{ \nti^2 } { ( \nti -1 )^2 } \frac{1}{ \ni - 1 }  \sumj  ( d_{ij} -  d_{i\cdot} ) ^2 \Big( \frac{1}{\nti} - \frac{1}{\ni} \Big)  \nonumber \\
& \leq &  4  \maxj  ( b_{ij} - b_{i\cdot} )^2 \frac{ \nci }{ \ni \nti } \S^2_{id}.
\end{eqnarray}
where the last inequality is due to $\nti^2 / (\nti - 1 )^2 \leq 4 $ ($\nti \geq 2$).

To bound the second term, we first observe that by Cauchy-Schwarz inequalty
\begin{eqnarray}
 ( \hat b_{i\cdot } - b_{i\cdot} )^2 = \Big\{    \frac{1}{\nti}  \sumj Z_{ij} ( b_{ij} - b_{i\cdot} )  \Big\}^2 \leq  \frac{1}{\nti} \sumj ( b_{ij} - b_{i\cdot} )^2 \leq \maxj ( b_{ij} - b_{i\cdot} )^2. \nonumber
\end{eqnarray}
Thus, the second term is bounded as follows:
\begin{eqnarray}
\label{eqn:second_term}
&&  \var \Big(  \frac{ \nti }{ \nti - 1} ( \hat b_{i\cdot } - b_{i\cdot} ) ( \hat d_{i \cdot} -  d_{i\cdot} )  \Big)  \nonumber \\
& \leq &  \frac{ \nti^2 } { ( \nti -1 )^2 } E \Big (  ( \hat b_{i\cdot } - b_{i\cdot} )^2 ( \hat d_{i \cdot} -  d_{i\cdot} )^2   \Big) \nonumber \\
& \leq &  \maxj  ( b_{ij} - b_{i\cdot} )^2   \frac{ \nti^2 } { ( \nti -1 )^2 }  E\Big( \frac{1}{ \nti }  \sumj  Z_{ij} ( d_{ij} -  d_{i\cdot} )   \Big)^2   \nonumber \\
& = & \maxj  ( b_{ij} - b_{i\cdot} )^2 \frac{ \nti^2 } { ( \nti -1 )^2 } \S^2_{id} \Big( \frac{1}{\nti} - \frac{1}{\ni} \Big)  \nonumber \\
& \leq &  4  \maxj  ( b_{ij} - b_{i\cdot} )^2 \frac{ \nci }{ \ni \nti } \S^2_{id}.
\end{eqnarray}
Combining \eqref{eqn:first_term} and \eqref{eqn:second_term}, we have
\[
\var( \s_{ibd}  )  \leq  16 \maxj  ( b_{ij} - b_{i\cdot} )^2 \frac{ \nci }{ \ni \nti } \S^2_{id}.
\]
Therefore,
\begin{eqnarray}
\var \Big( \sumi c_i \s_{ibd} \Big) & \leq & 16 \sumi c_i^2 \maxj  ( b_{ij} - b_{i\cdot} )^2 \frac{ \nci }{ \ni \nti } \S^2_{id} \nonumber \\
& \leq & 16 \frac{1}{N} \maxi \maxj   ( b_{ij} - b_{i\cdot} )^2 \sumi c_i \frac{ \nci }{ \nti } \S^2_{id},
\end{eqnarray}
which converges to zero by Condition \ref{cond:prop} and Condition \ref{cond:lem1}.
\end{proof}

\subsection*{Proof of lemma~\ref{lem:consistent_beta}}

\begin{proof}

By definition, 
\begin{eqnarray}
\hat \beta_1 & = &  \argmin\limits_{\beta}  \sumi \frac{ c_i }{n_{1i} - 1} \sumj Z_{ij}  \Big[  y_{ij}(1) - \hat y_{i\cdot}(1) - \big\{ \X_{ij} - \hat \X_{i\cdot}(1) \big\}^\T \beta \Big]^2. \nonumber \\
& = & \Big\{ \sumi c_i \s_{i\X\X}(1)  \Big\}^{-1} \Big\{ \sumi c_i \s_{i \X y}(1) \Big\},
\end{eqnarray}
where $ \s_{i\X\X}(1) $ is the sample covariance matrix of \revise{covariates} in stratum $i$, and $\s_{i \X y} (1)$ is the sample covariance between $\X_{ij}$ and $y_{ij}(1)$  in stratum $i$. Applying Lemma~\ref{lem:lem1} to $b_{ij} = \X_{ijk}$ and $d_{ij} = \X_{ijk'}$, we have $ \sumi c_i \s_{i\X\X}(1)  -  \sumi c_i \S_{i\X\X} $ and $\sumi c_i \s_{i \X y}(1)  - \sumi c_i \S_{i \X y}(1) $ converge, element-wise, to zero in probability. Moreover, by definite of $\beta_1$, we have
\[
\beta_1 =  \Big( \sumi c_i \S_{i\X\X}  \Big)^{-1} \Big( \sumi c_i \S_{i \X y}(1) \Big).
\]
Therefore, $\hat \beta_1 - \beta_1$ converges to zero in probability. Similarly, $\hat \beta_0 - \beta_0$ converges to zero in probability.

\end{proof}

\end{document}